\documentclass[12pt]{amsart} 
\usepackage[mathscr]{eucal}
\usepackage{amsmath,amsfonts}

\parskip=\smallskipamount

\newtheorem{theorem}{Theorem}[section]

\makeatletter
\@addtoreset{equation}{section}
\makeatother

\newcommand{\CC}{{\mathbb C}}
\newcommand{\NN}{{\mathbb N}}

\newcommand{\ZZ}{{\mathbb Z}}
\newcommand{\DD}{{\mathbb D}}
\newcommand{\RR}{{\mathbb R}}
\newcommand{\FF}{{\mathbb F}}
\newcommand{\TT}{{\mathbb T}}

\newcommand{\cA}{{\mathcal A}}
\newcommand{\cB}{{\mathcal B}}

\newcommand{\cD}{{\mathcal D}}
\newcommand{\cE}{{\mathcal E}}
\newcommand{\cF}{{\mathcal F}}

\newcommand{\cH}{{\mathcal H}}

\newcommand{\cL}{{\mathcal L}}
\newcommand{\cM}{{\mathcal M}}
\newcommand{\cN}{{\mathcal N}}
\newcommand{\cP}{{\mathcal P}}
\newcommand{\cR}{{\mathcal R}}
\newcommand{\cS}{{\mathcal S}}
\newcommand{\cU}{{\mathcal U}}

\newcommand{\cT}{{\mathcal T}}

\newdimen\expt
\def\boxit#1{\setbox0\hbox{$\displaystyle{#1}$}
      \hbox{\lower.4\expt
 \hbox{\lower3\expt\hbox{\lower\dp0
      \hbox{\vbox{\hrule height.4\expt
 \hbox{\vrule width.4\expt\hskip3\expt
      \vbox{\vskip3\expt\box0\vskip2\expt}%
 \hskip3\expt\vrule width.4\expt}\hrule height.4\expt}}}}}}
\begin{document}
\pagestyle{plain}

\bigskip

\title 
{Relations on noncommutative variables  \\
and associated orthogonal polynomials} 
\author{T. Banks} \author{T. Constantinescu} \author{J. L. Johnson} 

\address{Department of Mathematics \\
  University of Texas at Dallas \\
  Richardson, TX 75083}
\email{\tt banks@utdallas.edu} 
\address{Department of Mathematics \\
  University of Texas at Dallas \\
  Richardson, TX 75083} 
\email{\tt tiberiu@utdallas.edu}
\address{Department of Mathematics and Computer Science \\
  Wagner College \\
  Staten Island, NY 10301}
\email{\tt joeljohn@wagner.edu}

\noindent
\begin{abstract}
This semi-expository paper surveys results concerning
three classes of orthogonal polynomials: in one non-hermitian variable, 
in several isometric non-commuting variables, and in several 
hermitian non-commuting variables. The emphasis is on some dilation theoretic
techniques that are also described in some details.
\end{abstract}

\maketitle

\section{Introduction}
In this semi-expository paper we deal with a few classes of orthogonal 
polynomials 
associated to polynomial relations on several noncommuting variables.
Our initial interest in this subject was motivated by the need of more
examples related to \cite{CG}. 
On the other hand there are transparent connections with interpolation 
problems in several variables as well as with the modeling of various classes
of nonstationary systems (see \cite{AM} for a list of recent
references), which guided our choice of topics. 
Thus we do not relate to more traditional studies on orthogonal
polynomials of several variables associated to a finite reflection
group on an Euclidean space or other types of special functions of 
several variables, for which
a recent presentation could be found in \cite{DX}, instead we are more
focused on results connected with various dilation theoretic aspects
and Szeg\"o kernels. 

Our aim is to give an introduction to this point of view.
We begin our presentation with a familiar setting for algebras given 
by polynomial defining relations and then we introduce families of
orthonormal polynomials associated to some positive functionals on these
algebras. Section ~3 contains a discussion of the first class of 
orthogonal polynomials considered in this paper, namely  
polynomials in one variable on which there is no relation. This choice is 
motivated mainly by the fact that we have an opportunity to introduce 
some of the basic dilation theoretic techniques that we are using. First, 
we discuss (in a particular case that is sufficient for our goals) the 
structure of positive definite kernels and their triangular factorization.
Then these results are used to obtain recurrence relations for 
orthogonal polynomials in one variable with no additional relations, 
as well as asymptotic properties of these polynomials. All of these 
extend well-known results of Szeg\"o. 
We conclude this section with the introduction of a Szeg\"o type 
kernel which appears to be relevant to our setting.

In Section ~4 we discuss orthogonal polynomials of several
isometric variables. Most of the results are just 
particular cases of the corresponding results discussed in Section ~3, 
but there is an interesting new point about the Szeg\"o kernel that 
appears in the proof of Theorem ~4.1. We also use a certain explicit
structure of the Kolmogorov decomposition of a positive definite
kernel on the set of non-negative integers in order to
produce examples of families of operators satisfying 
Cuntz-Toeplitz and Cuntz relations.

The final section contains a discussion of orthogonal polynomials
of several non-commuting hermitian variables. This time, 
some of the techniques described in Section ~3 are not so relevant
and instead we obtain three-terms recursions in the traditional way, 
and we introduce families of Jacobi matrices associated to these 
recursions. 
Many of these results can be proved by adapting the classical proofs
from the one scalar variable case. However, much of the classical function 
theory is no longer available so we present some proofs illustrating how
classical techniques have to be changed or replaced. 
Also some results are not presented
in the most general form in the hope that the consequent simplifications 
in notation would make the paper more readable.

\section{Orthogonal polynomials associated to polynomial relations}
In this section we introduce some classes 
of orthogonal polynomials in several variables. We begin with the 
algebra $\cP _N$ of polynomials in $N$ noncommuting
variables  $X_1,\ldots ,X_N$ with complex coefficients.
Let $\FF _N^+$ be the unital free semigroup 
on $N$ generators $1,\ldots ,N$. The empty word is the identity element
and the length of the word $\sigma $ is 
denoted by $|\sigma |$. It is convenient to use the notation 
$X_{\sigma }=X_{i_1}\ldots X_{i_k}$ for $\sigma =i_1\ldots i_k
\in \FF _N^+$. Thus, each element $P\in \cP _N$ can be uniquely written 
in the form $P=\sum _{\sigma \in \FF _N^+}c_{\sigma }X_{\sigma }$,
with $c_{\sigma }\ne 0$ for finitely many $\sigma $'s.

We notice that $\cP _N$ is isomorphic with the tensor algebra 
over $\CC ^N$. Let $(\CC ^N)^{\otimes k}$ denote the $k$-fold tensor
product of $\CC ^N$ with itself. The tensor algebra 
over $\CC ^N$ is defined by the algebraic direct sum 
$$\cT (\CC ^N)=\oplus _{k\geq 0}(\CC ^N)^{\otimes k}.$$
If $\{e_1,\ldots ,e_N\}$ is the standard basis of
$\CC ^N$, then the set
$$\{1\}\cup \{\{e_{i_1}\otimes \ldots \otimes e_{i_k}  
\mid 1\leq i_1, \ldots ,i_k\leq N, k\geq 1\}
$$ 
is a  
basis of $\cT (\CC ^N)$. For $\sigma =i_1\ldots i_k$ we write
$e_{\sigma }$ instead of $e_{i_1}\otimes \ldots \otimes e_{i_k}$, 
and the mapping $X_{\sigma }\rightarrow e_{\sigma }$, $\sigma \in \FF _N^+$, 
extends to an isomorphism from  $\cP _N$ to $\cT (\CC ^N)$, 
hence $\cP _N \simeq \cT (\CC ^N)$

It is useful to introduce a natural involution on 
$\cP _{2N}$ as follows:
$$X^+_k=X_{N+k},\quad k=1,\ldots ,N,$$
$$X^+_l=X_{l-N},\quad l=N+1,\ldots ,2N;$$
on monomials, 
$$(X_{i_1}\ldots X_{i_k})^+=X^+_{i_k}\ldots X^+_{i_1},$$
and finally, if $Q=\sum _{\sigma \in \FF _{2N}^+}c_{\sigma }X_{\sigma }$,
then $Q^+=
\sum _{\sigma \in \FF _{2N}^+}\overline{c}_{\sigma }
X^+_{\sigma }$.
Thus, $\cP _{2N}$ is a unital, associative, $*$-algebra over $\CC $.

We say that $\cA \subset \cP _{2N}$
is  {\it symmetric} if $P\in \cA $ implies 
$cP^+\in \cA $ for some $c\in \CC -\{0\}$. 
Then the quotient of $\cP _{2N}$
by the two-sided ideal generated by 
$\cA $ is an associative algebra $\cR (\cA )$.
Letting $\pi =\pi _{\cA }:\cP _{2N}\rightarrow 
\cR (\cA )$ denote the quotient map then the formula 
\begin{equation}\label{invo}
\pi  (P)^+=\pi (P^+)
\end{equation}
gives a well-defined involution on $\cR (\cA )$.
Usually the elements of $\cA $ are called the {\it defining relations}
of the algebra $\cR (\cA )$. For instance, $\cR (\emptyset )=\cP _{2N}$,
$$\cR (\{X_k-X^+_k\mid k=1,\ldots ,N\})=\cP _N,$$
$$\cR (\{X_kX_l-X_lX_k\mid k,l=1,\ldots ,2N\})\simeq \cS (\CC ^{2N}),
\quad \mbox{the symmetric algebra over $\CC ^{2N}$,}$$
$$\cR (\{X_kX_l+X_lX_k\mid k,l=1,\ldots ,2N\})\simeq \Lambda (\CC ^{2N}),
\quad \mbox{the exterior algebra over $\CC ^{2N}$}.$$
Examples abound in the literature (for instance, see \cite{Dr}, \cite{DX},
\cite{Li}).

There are many well-known difficulties in the study of orthogonal
polynomials in several variables. The first one concerns the choice of an 
ordering of $\FF ^+_N$. In this paper we consider  
only the lexicographic
order $\prec $, but due to the canonical grading of 
$\FF _N^+$ it is possible to develop a basis free approach to 
orthogonal polynomials. In the case of orthogonal polynomials
on several commuting variables this is presented in \cite{DX}. 
A second difficulty concerns the choice of the moments. In this paper we 
adopt the following terminology. A linear functional $\phi $
on $\cR (\cA )$ is called {\it q-positive} (q comes from quarter) if 
$\phi (\pi (P)^+\pi (P))\geq 0$
for all $P\in \cP _N$. In this case, 
$\phi (\pi (P)^+)=
\overline{\phi (\pi (P))}$ for $P\in \cP _N$
and 
$$|\phi (\pi (P_1)^+\pi (P_2))|^2\leq 
\phi (\pi (P_1)^+\pi (P_1))
\phi (\pi (P_2)^+\pi (P_2))
$$
for $P_1,P_2\in \cP _N$.
Next we introduce 
\begin{equation}\label{scal}
\langle \pi (P_1),\pi (P_2)\rangle _{\phi }=
\phi (\pi (P_2)^+\pi (P_1)), \quad P_1,P_2\in \cP _N.
\end{equation}
By factoring out the subspace
$\cN _{\phi }=\{\pi (P)\mid P\in \cP _N,
\langle \pi (P),\pi (P) \rangle _{\phi }=0\}$
and 
completing this quotient with respect
to the norm induced by \eqref{scal} we obtain a Hilbert space 
$\cH _{\phi }$. 

The {\em index set} $G(\cA )\subset \FF ^+_N$ of $\cA$, 
is  
chosen as follows:
if $\alpha \in G(\cA )$, choose the next
element in $G(\cA )$ to be the least $\beta \in \FF _N^+$ 
with the property that the 
elements $\pi (X_{\alpha '})$, $\alpha '\preceq \alpha $, 
and $\pi (X_{\beta })$ are linearly independent.
We will avoid the degenerate situation in which $\pi (1)=0$;
if we do so, then $\emptyset \in G(\cA )$. 
Define $F_{\alpha }=\pi (X_{\alpha })$ for $\alpha \in G(\cA )$.
For instance, $G(\emptyset )=\FF ^+_N$, in which case $F_{\alpha }=
X_{\alpha }$, $\alpha \in \FF ^+_N$. Also,
$$G(\{X_kX_l-X_lX_k\mid k,l=1,\ldots ,2N\})=
\{i_1\ldots i_k\in \FF ^+_N
\mid i_1\leq \ldots \leq i_k, k\geq 0\},$$
and 
$$G(\{X_kX_l+X_lX_k\mid k,l=1,\ldots ,2N\})=
\{i_1\ldots i_k\in \FF ^+_N
\mid i_1<\ldots <i_k, 0\leq k\leq N\}$$
(we use the convention that for $k=0$, $i_1\ldots i_k$ is the empty word).

We consider the {\em moments} of $\phi $ to be the numbers
\begin{equation}\label{mome}
s_{\alpha ,\beta }=
\phi (F^+_{\alpha }F_{\beta })=
\langle F_{\beta },F_{\alpha }\rangle _{\phi }, 
\quad \alpha ,\beta \in G(\cA ).
\end{equation}
The {\em kernel of moments} is given by 
 $K_{\phi }(\alpha ,\beta )=
s_{\alpha ,\beta }$, $\alpha ,\beta \in G(\cA )$. 
We notice that $\phi $ is a q-positive functional on 
$\cR (\cA )$ if and only if $K_{\phi }$ 
is a positive definite kernel on $G(\cA )$.
However, $K_{\phi }$ does not determine $\phi $ uniquely. 
One typical situation when $K_{\phi }$ determines $\phi $ is  
$\{X_k-X^+_k\mid k=1,\ldots ,N\}\subset \cA $; a more general example 
is provided by the Wick polynomials,
$$X_kX^+_l-a_{k,l}\delta _{k,l}-\sum _{m,n=1}^Nc_{k,l}^{m,n}X^+_mX_n,
\quad k,l=1,\ldots ,N,$$
where $a_{k,l}$, $c_{k,l}^{m,n}$ are complex numbers and $\delta _{k,l}$
is the Kronecker symbol. 

The moment problem is trivial in this framework since it is obvious that 
the numbers $s_{\alpha ,\beta }$, $\alpha ,\beta \in G(\cA )$,
are the moments of a q-positive functional on $\cR (\cA )$ if and only if the 
kernel $K(\alpha ,\beta )=s_{\alpha ,\beta }$, $\alpha ,\beta \in G(\cA )$,
is positive definite. 

We now introduce 
orthogonal polynomials in $\cR (\cA )$. Assume that 
$\phi $ is strictly q-positive on $\cR (\cA )$, that is,
$\phi (\pi (P)^+\pi (P))>0$ for $\pi (P)\ne 0$. In this case 
$\cN _{\phi }=\{0\}$ and 
$\pi (\cP _N)$ can be viewed as a subspace of $\cH _{\phi }$.
Also, $\{F_{\alpha }\}_{\alpha \in G(\cA )}$ is a linearly independent
family in $\cH _{\phi }$ and 
the 
Gram-Schmidt procedure gives a family $\{\varphi _{\alpha  }\}
_{\alpha \in G(\cA )}$ of elements in 
$\pi (\cP _N)$ such that 

\begin{equation}\label{bond1}
\varphi _{\alpha  }=
\sum _{\beta \preceq \alpha }a_{\alpha ,\beta }F_{\beta },
\quad a_{\alpha ,\alpha }>0;
\end{equation}
\begin{equation}\label{bond2}
\langle \varphi _{\alpha }, \varphi _{\beta }\rangle _{\phi }=
\delta _{\alpha ,\beta },
\quad \alpha ,\beta \in G(\cA ).
\end{equation}

\noindent
The elements $\varphi _{\alpha  }$, $\alpha \in G(\cA )$, will be called
the {\em orthonormal polynomials} associated to $\phi $.
An explicit formula for the orthonormal polynomials  
can be obtained in the same manner as in 
the classical, one scalar variable case. Thus, set 
\begin{equation}\label{Deh}
D_{\alpha }=\det\left[s_{\alpha ',\beta '}\right]_
{\alpha ',\beta '\preceq \alpha }>0, \quad \alpha \in G(\cA ),
\end{equation}  
and from now on $\tau -1$ denotes the 
predecessor of $\tau $ with respect to the 
lexicographic order on $\FF _N^+$, while $\sigma +1$
denotes the successor of $\sigma $. 
We have:
$\varphi _{\emptyset }=s_{\emptyset ,\emptyset }^{-1/2}$ and for 
$\emptyset \prec \alpha $,
\begin{equation}\label{crop1}
\varphi _{\alpha }=\frac{1}{\sqrt{D_{\alpha -1}D_{\alpha }}}
{\det \left[\begin{array}{c}
\left[s_{\alpha ',\beta '}\right]_
{\alpha '\prec \alpha ;
\beta '\preceq \alpha } \\
  \\
\begin{array}{ccc}
F_{\emptyset } & \ldots & F_{\alpha }
\end{array}
\end{array}
\right]},
\end{equation}
with an obvious interpretation of the determinant.
In the following sections we will discuss in more details orthonormal
polynomials associated to some simple defining relations.

\section{No relation in one variable}
This simple case allows us to illustrate some general techniques that
can be used in the study of orthonormal polynomials. We have $\cA =\emptyset $
and $N=1$, so $\cR (\cA )=\cP _1$. The index set is $\NN _0$,
the set of nonnegative integers, and $F_n=X_1^n$, $n\in \NN _0$.
The moment kernel of a q-positive functional on 
$\cP _1$ is $K_{\phi }(n,m)=\phi ((X_1^n)^+X_1^m)$, $n,m\in \NN _0$,
and we notice that there is no restriction on $K_{\phi }$
other than being positive definite. We now discuss some tools that can be used 
in this situation.

\bigskip
\noindent
{\it 3.1. Positive definite kernels on $\NN _0$.} 
We discuss a certain structure (and parametrization) of positive definite 
kernels on $\NN _0$. The nature of this structure is revealed by looking at
the simplest examples. First, we consider a strictly positive matrix
$$S=\left[\begin{array}{cc}
1 & a \\
a & 1 
\end{array}
\right], \quad a\in \RR .$$ 
This matrix gives a new inner product on $\RR ^2$
by the formula 
$$\langle x,y\rangle _S=\langle Sx,y\rangle ,\quad x,y\in \RR ^2,$$
where $\langle \cdot ,\cdot \rangle $ denotes the Euclidean inner product on 
$\RR ^2$. Let $\{e_1,e_2\}$ be the standard basis of $\RR ^2$. By renorming 
$\RR ^2$ with $\langle \cdot ,\cdot \rangle _S$ the angle 
between $e_1$ and $e_2$ was modified to the new angle $\theta =
\theta (e_1,e_2)$ such that 
\begin{equation}\label{311}
\cos \theta (e_1,e_2)=\frac{\langle e_1,e_2\rangle _S}
{\|e_1\|_S\|e_2\|_S}=a.
\end{equation}
We can visualize the renormalization process by giving a map
$T_S:\RR ^2\rightarrow \RR ^2$ with the property that 
$\langle T_Sx,T_Sy\rangle =\langle x,y\rangle _S$
for $x,y\in \RR ^2$,
and it is easily seen that we can choose 
$$T_S=
\left[\begin{array}{cc}
1 & \cos \theta  \\
0 & \sin \theta  
\end{array}
\right].$$

\begin{figure}[h]
\setlength{\unitlength}{2700sp}%
\begingroup\makeatletter\ifx\SetFigFont\undefined%
\gdef\SetFigFont#1#2#3#4#5{%
  \reset@font\fontsize{#1}{#2pt}%
  \fontfamily{#3}\fontseries{#4}\fontshape{#5}%
  \selectfont}%
\fi\endgroup%
\begin{picture}(4224,1824)(289,-1273)
{ \thinlines

{ \put(1400,-600){$e_1$}}
{ \put(900,-100){$e_2$}}
{ \put(3800,-600){$e_1$}}
{ \put(3950,180){$f_2$}}
{ \put(2400,350){$T_S$}}
{ \put(3800,-250){$\theta $}}

{\put(301,-361){\line( 1, 0){1800}}
}%
{\put(1201,539){\line( 0,-1){1800}}
}%
{\put(2701,-361){\line( 1, 0){1800}}
}%
{\put(3601,539){\line( 0,-1){1800}}
}%
{\put(1201,-361){\vector( 1, 0){600}}
}%
{\put(1201,-361){\vector( 0, 1){600}}
}%
{\put(3601,-361){\vector( 1, 0){600}}
}%
{\put(3601,-361){\vector( 1, 2){270}}
}%
\thicklines
{\multiput(1951,239)(184,0){7}{\line( 1, 0){ 92}}
\put(3151,239){\vector( 1, 0){0}}
}%
}
\end{picture}

\caption{\mbox{ Renormalization in $\RR ^2$ }}
\end{figure}

We can also notice that $T_Se_1=e_1$ and $T_Se_2=f_2=J(\cos \theta )e_1$, 
where $J(\cos \theta )$ is the Julia operator,
$$J(\cos \theta )=\left[\begin{array}{cc}
\cos \theta & \sin \theta \\
\sin \theta & -\cos \theta 
\end{array}\right],
$$
which is the composition of a reflection about the $x$-axis followed by
the counterclockwise rotation $R_{\theta }$ through angle $\theta $.
We deduce that 
$$a=\cos \theta =\langle e_1,f_2\rangle=\langle e_1,J(\cos \theta )e_1\rangle
=\langle J(\cos \theta )e_1,e_1\rangle .
$$

The discussion extends naturally to the $3\times 3$ case. Thus let
$$S=\left[\begin{array}{ccc}
1 & a & b \\
a & 1 & c \\
b & c & 1 
\end{array}\right], \quad a,b,c\in \RR,
$$
be a strictly positive matrix. A new inner product is induced by 
$S$ on $\RR ^3$, 
$$\langle x,y\rangle _S=\langle Sx,y\rangle ,\quad x,y\in \RR ^3,$$
and let $\{e_1,e_2,e_3\}$ be the standard basis of $\RR ^3$. With respect to 
this new
inner product the vectors $e_1$, $e_2$, $e_3$ still belong to the unit
sphere, but they are no longer orthogonal. Thus, 
$$a=\cos \theta (e_1,e_2)=\cos \theta _{12},$$
$$c=\cos \theta (e_2,e_3)=\cos \theta _{23},$$
and
$$b=\cos \theta (e_1,e_3)=\cos \theta _{13}.$$
This time, the law of cosines in spherical geometry gives a relation
between the numbers $a$, $b$, and $c$, 
\begin{equation}\label{cosi}
b=\cos \theta _{13}=\cos \theta _{12}\cos \theta _{23}+
\sin \theta _{12}\sin \theta _{23}\cos \theta ,
\end{equation}
where $\theta $ is the dihedral angle formed by the planes
generated by $e_1$, $e_2$ and, respectively, $e_2$, $e_3$
(see, for instance, \cite{HH}). Thus, the number $b$ belongs to a disk of 
center $\cos \theta _{12}\cos \theta _{23}$ and radious 
$\sin \theta _{12}\sin \theta _{23}$.
Once again the renormalization can
be visualized by a map $T_S:\RR ^3\rightarrow \RR ^3$
such that $\langle T_Sx,T_Sy\rangle =\langle x,y\rangle _S$. In this case we 
can choose
$$T_S=\left[\begin{array}{ccc}
1 & \cos \theta _{12} & \cos \theta _{12}\cos \theta _{23}+
\sin \theta _{12}\sin \theta _{23}\cos \theta \\
0 & \sin \theta _{12} & \sin \theta _{12}\cos \theta _{23}-
\cos \theta _{12}\sin \theta _{23}\cos \theta \\ 
0 & 0 & \sin \theta _{23}\sin \theta
\end{array}\right],$$
and we see that $T_Se_1=e_1$, $T_Se_2=f_2=(J(\cos \theta _{1,2})\oplus 1)e_1$, 
and
$$T_Se_3=f_3=
(J(\cos \theta _{12})\oplus 1)
(1\oplus J(\cos \theta ))
(J(\cos \theta _{23})\oplus 1)e_1.
$$ 
In particular, 
\begin{equation}\label{dila}
b=\cos \theta _{13}=
\langle (J(\cos \theta _{1,2})\oplus 1)
(1\oplus J(\cos \theta ))
(J(\cos \theta _{2,3})\oplus 1)e_1,e_1\rangle ,
\end{equation}
which can be viewed as a dilation formula.

Now both \eqref{cosi} and \eqref{dila} extend to a 
strictly q-positive $n\times n$ matrix and provide 
a parametrization and therefore a structure for the positive definite
kernels on $\NN _0$ (for general results and details see \cite{Co2}, 
\cite{Co}). We apply this result to a kernel $K_{\phi }$ associated 
to a strictly q-positive functional $\phi $ and obtain that 
$K_{\phi }$ is 
uniquely determined by a family $\{\gamma _{k,j}\}_{0\leq k<j}$
of complex numbers with the property that $|\gamma _{k,j}|<1$
for all $0\leq k<j$.  
Define $d_{k,j}=(1-|\gamma _{k,j}|^2)^{1/2}$. The extension of 
\eqref{dila} mentioned above gives
\begin{equation}\label{gdil}
s_{k,j}=s^{1/2}_{k,k}s^{1/2}_{j,j}\langle U_{k,j}e_1,e_1\rangle ,\quad k<j,
\end{equation}
where $U_{k,j}$ is a $(j-k+1)\times (j-k+1)$ 
unitary matrix defined recursively by: $U_{k,k}=1$, 
and for $k<j$, 
\begin{equation}\label{us}
U_{k,j}=
(J(\gamma _{k,k+1})\oplus 1_{n-1})
(1\oplus J(\gamma _{k,k+2})\oplus 1_{n-2})\ldots 
(1_{n-1}\oplus J(\gamma _{k,j}))
(U_{k+1,j}\oplus 1).
\end{equation}
Also $J(\gamma _{l,m})$ is the Julia operator associated to 
$\gamma _{l,m}$, that is, 
$$J(\gamma _{l,m})=\left[\begin{array}{cc}
\gamma _{l,m} & d_{l,m} \\
d_{l,m} & -\overline{\gamma }_{l,m}
\end{array}\right],$$
and $1_m$ denotes the $m\times m$ identity matrix.
For instance, we deduce from \eqref{gdil} that:
$$s_{01}=s^{1/2}_{00}s^{1/2}_{11}\left[
\begin{array}{cc}
1 & 0   
\end{array}
\right]
\left[
\begin{array}{cc}
\gamma _{01} & d_{01} \\
d_{01} & -\overline{\gamma }_{01} 
\end{array}
\right]
\left[
\begin{array}{c}
1 \\
 0 
\end{array}
\right];
$$
$$
\begin{array}{rcl}
s_{02}\!\!\!&=&\!\!\!s^{1/2}_{00}s^{1/2}_{22}\left[
\begin{array}{ccc}
1 & 0 & 0  
\end{array}
\right] \\
 & & \\
 & & \quad \quad \times \left[
\begin{array}{ccc}
\gamma _{01} & d_{01} & 0 \\
d_{01} & -\overline{\gamma }_{01} & 0 \\
 0 & 0 & 1
\end{array}
\right]
\left[
\begin{array}{ccc}
 1 & 0 & 0 \\
 0 & \gamma _{02} & d_{02} \\
 0 & d_{02} & -\overline{\gamma }_{02} 
\end{array}
\right]
\left[
\begin{array}{ccc}
\gamma _{12} & d_{12} & 0 \\
d_{12} & -\overline{\gamma }_{12} & 0 \\
 0 & 0 & 1
\end{array}
\right]
\left[
\begin{array}{c}
1 \\
 0 \\
 0 
\end{array}
\right],
\end{array}
$$
and 
$$\begin{array}{rcl}
s_{03}\!\!\!&=&\!\!\!s^{1/2}_{00}s^{1/2}_{33}\left[
\begin{array}{cccc}
1 & 0 & 0 & 0 
\end{array}
\right] \\
 & & \\
 & & \quad \quad \times \left[
\begin{array}{cccc}
\gamma _{01} & d_{01} & 0 & 0 \\
d_{01} & -\overline{\gamma }_{01} & 0 & 0 \\
 0 & 0 & 1 & 0 \\
 0 & 0 & 0 & 1
\end{array}
\right]
\left[
\begin{array}{cccc}
 1 & 0 & 0 & 0 \\
 0 & \gamma _{02} & d_{02} & 0 \\
 0 & d_{02} & -\overline{\gamma }_{02} & 0 \\
  0 & 0 & 0 & 1
\end{array}
\right]
\left[
\begin{array}{cccc}
 1 & 0 & 0 & 0 \\
 0 & 1 & 0 & 0 \\
 0 & 0 & \gamma _{03} & d_{03} \\
 0 & 0 & d_{03} & -\overline{\gamma }_{03} \\
\end{array}
\right] \\
 & & \\
 & &  \quad \quad \times \left[ 
\begin{array}{cccc}
\gamma _{12} & d_{12} & 0 & 0 \\
d_{12} & -\overline{\gamma }_{12} & 0 & 0 \\
 0 & 0 & 1 & 0 \\
 0 & 0 & 0 & 1
\end{array}
\right]
\left[
\begin{array}{cccc}
 1 & 0 & 0 & 0 \\
 0 & \gamma _{13} & d_{13} & 0 \\
 0 & d_{13} & -\overline{\gamma }_{13} & 0 \\
  0 & 0 & 0 & 1
\end{array}
\right] \\
 & & \\
 & & \quad \quad \times \left[
\begin{array}{cccc}
\gamma _{23} & d_{23} & 0 & 0 \\
d_{23} & -\overline{\gamma }_{23} & 0 & 0 \\
 0 & 0 & 1 & 0 \\
 0 & 0 & 0 & 1
\end{array}
\right]
\left[
\begin{array}{c}
1 \\
 0 \\
 0 \\
 0
\end{array}
\right]. 
\end{array}
$$
In particular, we deduce that 
$s_{01}=s^{1/2}_{00}s^{1/2}_{11}\gamma _{01}$. The next formula is 
the complex version of \eqref{cosi}, 
$s_{02}=s^{1/2}_{00}s^{1/2}_{22}\left(
\gamma _{01}\gamma _{12}+d_{01}\gamma _{02}d_{12}\right)$.
Then, 
$$\begin{array}{rcl}
s_{03}\!\!\!&=&\!\!\!s^{1/2}_{00}s^{1/2}_{33}
\left(\gamma _{01}\gamma _{12}\gamma _{23}+
\gamma _{01}d_{12}\gamma _{13}d_{23}
+d_{01}\gamma _{02}d_{12}\gamma _{23}\right.\\
 & & \quad \left. 
-d_{01}\gamma _{02}\overline{\gamma }_{12}\gamma _{13}d_{23}+
d_{01}d_{02}\gamma _{03}d_{13}d_{23}\right).
\end{array}
$$
Explicit formulae of this type can be obtained for each $s_{k,j}$,
as well as inverse algorithms allowing to calculate $\gamma _{k,j}$
from the kernel of moments, see \cite{Co} for details.

A natural combinatorial question would be to calculate the number 
$N(s_{k,j})$ of additive terms in the expression of $s_{k,j}$.
We give here some details since the calculation of $N(s_{k,j})$ involves 
another useful interpretation of formula \eqref{gdil}.   
We notice that for $k\geq 0$,
$$N(s_{01})=N({s_{k,k+1}})=1,$$
$$N(s_{02})=N({s_{k,k+2}})=2,$$
$$N(s_{03})=N({s_{k,k+3}})=5.$$
The general formula is given by the following result.
\begin{theorem}\label{catalan}
$N(s_{k,k+l})$ is given by the Catalan number $C_l=\displaystyle\frac{1}{l+1}
\left(\begin{array}{c}
2l \\
l 
\end{array}
\right)$.
\end{theorem}
\begin{proof}
The first step of the proof considers the realization of $s_{k,j}$
through a time varying transmission line (or lattice). 
For illustration we consider the case of $s_{03}$
in Figure ~2.

\begin{figure}[h]
\setlength{\unitlength}{2700sp}%
\begingroup\makeatletter\ifx\SetFigFont\undefined%
\gdef\SetFigFont#1#2#3#4#5{%
  \reset@font\fontsize{#1}{#2pt}%
  \fontfamily{#3}\fontseries{#4}\fontshape{#5}%
  \selectfont}%
\fi\endgroup%
\begin{picture}(6324,2874)(289,-2323)
{ \thinlines
}%
{ \put(301,-361){\line( 1, 0){900}}
}%
{ \put(1201,-361){\vector( 1,-1){600}}
}%
{ \put(1201,-961){\vector( 1, 1){600}}
}%
{ \put(500,10){$A$}}
{ \put(6000,10){$B$}}
{ \put(1126,-1111){\framebox(750,900){}}
}%
{ \put(901,-961){\line( 1, 0){1200}}
}%
{ \put(1201,-361){\line( 1, 0){1200}}
}%
{ \put(2926,-1111){\framebox(750,900){}}
}%
{ \put(2026,-1711){\framebox(750,900){}}
}%
{ \put(3826,-1711){\framebox(750,900){}}
}%
{ \put(4726,-1111){\framebox(750,900){}}
}%
{ \put(2926,-2311){\framebox(750,900){}}
}%
{ \put(2101,-961){\line( 1, 0){1200}}
}%
{ \put(3301,-961){\line( 1, 0){1200}}
}%
{ \put(4501,-961){\line( 1, 0){1200}}
}%
{ \put(2401,-361){\line( 1, 0){1200}}
}%
{ \put(3601,-361){\line( 1, 0){1200}}
}%
{ \put(4801,-361){\line( 1, 0){1200}}
}%
{ \put(1201,-361){\vector( 1, 0){600}}
}%
{ \put(1201,-961){\vector( 1, 0){600}}
}%
{ \put(3001,-361){\vector( 1, 0){600}}
}%
{ \put(3001,-961){\vector( 1, 0){600}}
}%
{ \put(2101,-961){\vector( 1, 0){600}}
}%
{ \put(2101,-1561){\vector( 1, 0){600}}
}%
{ \put(3001,-1561){\vector( 1, 0){600}}
}%
{ \put(3001,-2161){\vector( 1, 0){600}}
}%
{ \put(3901,-961){\vector( 1, 0){600}}
}%
{ \put(3901,-1561){\vector( 1, 0){600}}
}%
{ \put(4801,-361){\vector( 1, 0){600}}
}%
{ \put(4801,-961){\vector( 1, 0){600}}
}%
{ \put(3001,-961){\vector( 1, 1){600}}
}%
{ \put(4801,-961){\vector( 1, 1){600}}
}%
{ \put(3001,-361){\vector( 1,-1){600}}
}%
{ \put(4801,-361){\vector( 1,-1){600}}
}%
{ \put(2101,-1561){\vector( 1, 1){600}}
}%
{ \put(3901,-1561){\vector( 1, 1){600}}
}%
{ \put(3001,-2161){\vector( 1, 1){600}}
}%
{ \put(2176,-961){\vector( 1,-1){600}}
}%
{ \put(3001,-1561){\vector( 1,-1){600}}
}%
{ \put(3901,-961){\vector( 1,-1){600}}
}%
{ \put(1801,-1561){\line( 1, 0){3000}}
}%
{ \put(2701,-2161){\line( 1, 0){1200}}
}%
\end{picture}

\caption{\mbox{ Lattice representation for $s_{03}$}}
\end{figure}

Each box in Figure~2 represents the action of a Julia operator,
and we see on this figure that the number of additive terms in the formula of 
$s_{03}$ is given by
the number of paths from $A$ to $B$. In it clear that 
to each 
path from $A$ to $B$ in Figure ~2 it corresponds a Catalan path from 
$C$ to $D$ in Figure ~3, that is, a path that never steps 
below the diagonal and goes only to the right or 
downward. 

\begin{figure}[h]
\setlength{\unitlength}{2700sp}%
\begingroup\makeatletter\ifx\SetFigFont\undefined%
\gdef\SetFigFont#1#2#3#4#5{%
  \reset@font\fontsize{#1}{#2pt}%
  \fontfamily{#3}\fontseries{#4}\fontshape{#5}%
  \selectfont}%
\fi\endgroup%
\begin{picture}(1800,1733)(76,-1111)
{\thinlines
\put(751,539){\circle{150}}
}%
{ \put(1201,539){\circle{150}}
}%
{\put(-100,539){{$C$}}
}%
{\put(1900,-811){{$D$}}
}%
{ \put(1651,539){\circle{150}}
}%
{ \put(301, 89){\circle{150}}
}%
{ \put(751, 89){\circle{150}}
}%
{ \put(1201, 89){\circle{150}}
}%
{ \put(1651, 89){\circle{150}}
}%
{ \put(301,-361){\circle{150}}
}%
{ \put(751,-361){\circle{150}}
}%
{ \put(1201,-361){\circle{150}}
}%
{ \put(1651,-361){\circle{150}}
}%
{ \put(301,-811){\circle{150}}
}%
{ \put(751,-811){\circle{150}}
}%
{ \put(1201,-811){\circle{150}}
}%
{ \put(1651,-811){\circle{150}}
}%
{ \put(301,539){\circle{150}}
}%
{ \put(301,539){\line( 1, 0){450}}
\put(751,539){\line( 0,-1){450}}
\put(751, 89){\line( 1, 0){900}}
\put(1651, 89){\line( 0,-1){900}}
}%
\end{picture}

\caption{\mbox{ A Catalan path from $C$ to $D$}}
\end{figure}

Thus, each box in Figure~2 corresponds to a point 
strictly above the diagonal 
in Figure~3. Once this one-to-one correspondence is established, 
we can use the 
well-known fact that the number of Catalan paths like the one in Figure~3 
is given 
exactly by the Catalan numbers.
\end{proof}

\bigskip
\noindent
{\it 3.2. Spectral factorization.}
The classical theory of orthogonal polynomials is intimately related
to the so-called spectral factorization. Its prototype would be the 
Fej\'er-Riesz factorization of a positive trigonometric 
polynomial $P$ in the form $P=|Q|^2$, where $Q$ is a 
polynomial with no zeros in the unit disk. This is 
generalized by Szeg\"o to the Szeg\"o class of those measures on 
the unit circle $\TT $ with $\log \mu '\in L^1$, and very
general results along this line can be found in \cite{RR} and \cite{NF}. 

Here we briefly review the spectral factorization of positive definite 
kernels on the set $\NN _0$ described in \cite{Co4}.
For two positive definite kernels $K_1$ and $K_2$
we write $K_1\leq K_2$ if $K_2-K_1$
is a positive definite kernel.
Consider a family 
$\cF =\{\cF _n\}_{n\geq 0}$ of at most one-dimensional vector spaces and call 
{\sl lower triangular array}
a family $\Theta =\{\Theta _{k,j}\}_{k,j\geq 0}$ of complex numbers 
$\Theta _{k,j}$ with the following two properties: $\Theta _{k,j}=0$
for $k<j$ and each column
$c_j(\Theta )=[\Theta _{k,j}]_{k\geq 0}$, $j\geq 0$,  belongs
to the Hilbert space $\oplus _{k\geq j}\cF _k$. 
Denote by $\cH_0^2(\cF)$ the set of all lower triangular 
arrays as above. An element of  $\cH_0^2(\cF)$ 
is called {\sl outer} if the set  
$\{c_j(\Theta )\mid j\geq k\}$ is total in $\oplus _{j\geq k}\cF _j$
for each $k\geq 0$. 

It is easily seen that if $\Theta $ is an outer triangular array, then 
the formula
$$K_\Theta (k,j)=c_k(\Theta )^*c_j(\Theta )$$
gives a positive definite kernel on $\NN _0$.
The following result extends the above mentioned Szeg\"o factorization
and at the same time it contains the Cholesky factorization of positive 
matrices. 
\begin{theorem}\label{factorizare}
Let $K$ be an positive definite kernel on $\NN _0$. Then there exists a family
$\cF =\{\cF _n\}_{n\geq 0}$ of at most one-dimensional vector spaces
and an outer triangular array $\Theta \in \cH_0^2(\cF )$ such that 

\smallskip
\quad $(1)$ \quad $K_\Theta \leq K$.

\smallskip
\quad $(2)$ \quad For any other family 
$\cF '=\{\cF '_n\}_{n\geq 0}$ of at most one-dimensional vector spaces
and any outer triangular array $\Theta '\in \cH_0^2(\cE ,\cF ')$ such that 
$K_{\Theta '}\leq K$, we have $K_{\Theta '}\leq K_{\Theta }$. 

\smallskip
\quad $(3)$ \quad $\Theta $ is uniquely determined by $(a)$ and $(b)$ up to 
a left 
unitary diagonal factor.
\end{theorem}

It follows from $(3)$ above that the spectral factor $\Theta $ can be 
uniquely determined 
by the condition that $\Theta _{n,n}\geq 0$ for all $n\geq 0$.
We say that the kernel $K$ belongs to the Szeg\"o class if  
$\inf _{n\geq 0}\Theta _{n,n}>0$. If $\{\gamma _{k,j}\}$ are the 
parameters of $K$ introduced in Subsection ~3.1 then it follows that 
the kernel $K$ belongs to the Szeg\"o class if and only if 
\begin{equation}\label{36}
\inf _{k\geq 0}s_{k,k}^{1/2}\prod _{n>k}d_{k,n}>0.
\end{equation}
This implies that 
$\cF _n=\CC $ for all $n\geq 0$ (for details see \cite{Co4} or \cite{Co}).

\bigskip
\noindent
{\it 3.3. Recurrence relations.}
Formula \eqref{crop1} is not very useful in calculations 
involving the orthogonal polynomials. Instead there are used
recurrence formulae. In our case, $\cA =\emptyset $ and $N=1$, 
we consider the moment kernel $K_{\phi }$ of a strictly q-positive 
functional on $\cP _1$ and also, the parameters $\{\gamma _{k,j}\}$ 
of $K_{\phi }$ as in Subsection ~3.1.
It can be shown that the orthonormal polynomials associated to 
$\phi $ obey the following recurrence 
relations 
\begin{equation}\label{lazero}
\varphi _0(X_1,l)=\varphi _0^{\sharp }(X_1,l)=s_{l,l}^{-1/2}, 
\quad l\in \NN _0,
\end{equation}
and for $n\geq 1$, $l\in \NN _0$,
\begin{equation}\label{primarelatie}
\varphi _n(X_1,l)=\frac{1}{d_{l,n+l}}
\left( X_1\varphi _{n-1}(X_1,l+1)-
\gamma _{l,n+l}\varphi ^{\sharp }_{n-1}(X_1,l)\right),
\end{equation}
\begin{equation}\label{adouarelatie}
\varphi ^{\sharp}_n(X_1,l)=\frac{1}{d_{l,n+l}}
\left(-\overline{\gamma }_{l,n+l}X_1\varphi _{n-1}(X_1,l+1)+
\varphi ^{\sharp }_{n-1}(X_1,l)\right),
\end{equation}
where $\varphi _n(X_1)=\varphi _n(X_1,0)$ and 
$\varphi ^{\sharp }_n(X_1)=\varphi ^{\sharp }_n(X_1,0)$.

Somewhat similar polynomials are considered in \cite{DGK}, but the 
form of the recurrence relations as above is noticed in \cite{CJ}.
It should be mentioned that $\{\varphi _n(X_1,l)\}_{n\geq 0}$ is the 
family of orthonormal polynomials
associated to a q-positive functional on $\cP _1$ with moment
kernel $K^l(\alpha ,\beta )=
s_{\alpha +l,\beta +l}$,
$\alpha , \beta \in \NN _0$. 
Also, the above
recurrence 
relations provide us with a tool to recover the parameters 
$\{\gamma _{k,j}\}$ from the orthonormal polynomials.
\begin{theorem}\label{nouaparam}
Let $k^l_n$ be the leading coefficient of $\varphi _n(X_1,l)$.
For $l\in \NN _0$ and $n\geq 1$, 
$$\gamma _{l,n+l}=-
\varphi _n(0,l)\displaystyle\frac{k^{l+1}_0\ldots k^{l+1}_{n-1}}
{k^l_0\ldots k^l_n}.
$$
\end{theorem}
\begin{proof}
We reproduce here the proof from \cite{BC} in order to illustrate 
these concepts and to introduce one more property of the 
parameters $\{\gamma _{k,j}\}$.
First, we deduce
from \eqref{primarelatie} that
$$\varphi _n(0,l)=-\frac{\gamma _{l,n+l}}{d_{l,n+l}}
\varphi ^{\sharp}_{n-1}(0,l),
$$
while formula \eqref{adouarelatie} 
gives 
$$\varphi ^{\sharp }_n(0,l)=
\frac{1}{d_{l,n+l}}\varphi ^{\sharp }_{n-1}(0,l)=\ldots 
=s_{l,l}^{-1/2}\prod _{p=1}^n\frac{1}{d_{l,p+l}},
$$
hence
$$\varphi _n(0,l)=-s_{l,l}^{-1/2}\gamma _{l,n+l}
\prod _{p=1}^n\frac{1}{d_{l,p+l}}.
$$
Now we can use another useful feature of the parameters $\{\gamma _{k,j}\}$, 
namely the fact that they give simple formulae for determinants.
Let $D_{m,l}$ denote the determinant of the matrix 
$\left[s_{k,j}\right]_{l\leq k,j\leq m}$. By Proposition ~1.7 in \cite{Co2},
\begin{equation}\label{deter}
D_{l,m}=\prod _{k=l}^ms_{k,k}\times \prod _{l\leq j<p\leq m}d^2_{j,p}.
\end{equation}
One simple application of this formula is that it reveals the 
equality behind Fisher-Hadamard inequality. Thus, 
for $l\leq n\leq n'\leq m$, we have
$$D_{l,m}=\displaystyle\frac{D_{l,n'}D_{n,m}}
{D_{n',n}}\prod _{(k,j)\in \Lambda }d^2_{k,j},
$$
where $\Lambda =\{(k,j)\mid l\leq k<n\leq n'<j\leq m\}$. Some 
other applications of \eqref{deter} can be found in \cite{Co}, Chapter ~8.
Returning to our proof we deduce from \eqref{deter} that
$$\prod _{p=1}^nd^2_{l,p+l}=s^{-1}_{l,l}
\displaystyle\frac{D_{l,l+n}}{D_{l+1,l+n}}
$$
so, 
\begin{equation}\label{gama}
\gamma _{l,n+l}=-
\varphi _n(0,l)\sqrt{\frac{D_{l,l+n}}{D_{l+1,l+n}}}.
\end{equation}
We can now relate this formula to the leading coefficients $k^l_n$.
From \eqref{primarelatie} we deduce that
$$k^l_n=s^{-1/2}_{l+n,l+n}\prod _{p=1}^{n-1}\frac{1}{d_{l+p,l+n}},
\quad n\geq 1,$$
and using once again \eqref{deter}, we deduce
$$k^l_n=\sqrt{\frac{D_{l,l+n-1}}{D_{l,l+n}}},\quad n\geq 1,
$$
which concludes the proof.
\end{proof}

\bigskip
\noindent
{\it 3.4. Some examples.}
We consider some examples, especially in order to 
clarify the connection with classical orthogonal polynomials. Thus, 
consider first $\cA =\{1-X^+_1X_1\}$. 
In this case the index set is still $\NN _0$ 
and if $\phi $ is a linear functional 
on $\cR (\cA )$, then the kernel of moments is Toeplitz,  
$K_{\phi }(n+k,m+k)=K_{\phi }(n,m)$, $m,n,k\in \NN _0.$
Let $\phi $ be a strictly q-positive functional on $\cR (\cA )$
and let $\{\gamma _{k,j}\}$ be the parameters associated to $K_{\phi }$.
We deduce that these parameters also satisfy the Toeplitz condition, 
$\gamma _{n+k,m+k}=\gamma _{n,m}$, $n<m$, $k\geq 1$. 
Setting $\gamma _n=\gamma _{k,n+k}$, $n\geq 1$, $k\geq 0$, 
and $d_n=(1-|\gamma _n|^2)^{1/2}$, the recurrence relations 
\eqref{primarelatie}, \eqref{adouarelatie} collaps to the classical
Szeg\"o recursions obeyed by the orthogonal polynomials on the unit
circle, 
$$
\varphi _{n+1}(z)=\frac{1}{d_{n+1}}(z\varphi _n(z)-
\gamma _{n+1}\varphi _n^{\sharp }(z)),
$$
and
$$
\varphi _{n+1}^{\sharp }(z)=\frac{1}{d_{n+1}}
(-\overline{\gamma }_{n+1}z\varphi _n(z)+
\varphi _n^{\sharp }(z)).
$$
Therefore $\gamma _n$, $n\geq 1$ are the usual Szeg\"o coefficients, 
\cite{Sz}.

Another example is given by
$\cA =\{X_1-X^+_1\}$. In this case the index set is still   
$\NN _0$ and the moment kernel of a strictly q-positive functional
on $\cR (\cA )$ will have the Hankel property, 
$K_{\phi }(n,m+k)=K_{\phi }(n+k,m)$, $m,n,k\in \NN _0.$
Orthogonal polynomials associated to functionals on $\cR (\cA )$ correspond
to orthogonal polynomials on the real line. This time, the parameters
$\{\gamma _{k,j}\}$ associated to moment kernels have no classical
analogue. Instead there are so-called canonical moments which are
used as a counterpart of the Szeg\"o coefficients (see \cite{He}).
Also, recurrence relations of type 
\eqref{primarelatie}, \eqref{adouarelatie} are replaced by a three term 
recurrence equation,
\begin{equation}\label{treitermeni}
x\varphi _n(x)=b_{n}\varphi _{n+1}(x)+
a_{n}\varphi _n(x)+
b_{n-1}\varphi _{n-1}(x),
\end{equation}
with initial conditions $\varphi _{-1}=0$, $\varphi _0=1$ (\cite{Sz}).
Definitely, these objects are more useful (for instance, it appears that
no simple characterization of those 
$\{\gamma _{k,j}\}$ corresponding to Hankel kernels is known). Still, 
computations involving the parameters $\{\gamma _{k,j}\}$ might be 
of interest. For instance, we show here how to calculate the
parameters for Gegenbauer polynomials. For a number  
$\lambda >-\frac{1}{2}$, these are orthogonal polynomials associated to
the weight function 
$w(x)=B(\frac{1}{2}, \lambda +\frac{1}{2})^{-1}(1-x^2)^{\lambda -\frac{1}{2}}
$
on $(-1,1)$ ($B$ denotes the beta function). 
We use the normalization constants from \cite{DX}, 
thus the Gegenbauer polynomials are
$$P^{\lambda }_n(x)=\frac{(-1)^n}{2^n(\lambda +\frac{1}{2})_n}
(1-x^2)^{\frac{1}{2}-\lambda }
\frac{d^n}{dx^n}(1-x^2)^{n+\lambda -\frac{1}{2}},$$
where $(x)_n$ is the Pochhammer symbol, 
$(x)_0=1$ and $(x)_n=\prod _{k=1}^n(x+k-1)$ for $n\geq 1$.
We have:
$$h^{\lambda }_n=\frac{1}{B(\frac{1}{2}, \lambda +\frac{1}{2})}
\int _{-1}^1\left(P^{\lambda }_n(x)\right)^2(1-x^2)^{\lambda -\frac{1}{2}}dx
=\frac{n!(n+2\lambda )}{2(2\lambda +1)_n(n+\lambda )}$$
and the three term recurrence is:
$$P^{\lambda }_{n+1}(x)=\frac{2(n+\lambda )}{n+2\lambda }xP^{\lambda }_n(x)
-\frac{n}{n+2\lambda }P^{\lambda }_{n-1}(x)
$$
(see \cite{DX}, Ch. 1). We now let $\varphi ^{\lambda }_n(x,0)$ denote the 
orthonormal polynomials associated to the weight function $w$,
hence $\varphi ^{\lambda }_n(x,0)=\frac{1}{\sqrt{h^{\lambda }_n}}
P^{\lambda }_n(x)$.
From the three term relation we deduce
$$\varphi ^{\lambda }_n(0,0)=(-1)^{n+1}
\sqrt{\frac{2(2\lambda +1)_n(n+\lambda )}{
n!(n+2\lambda )}}\times \prod _{k=1}^n\frac{k-1}{k-1+2\lambda },
$$
and also, the leading coefficient of $\varphi ^{\lambda }_n(x,0)$
is 
$$k^{\lambda ,0}_n=
\frac{(n+2\lambda )_n}{2^n\left(\lambda +\frac{1}{2}\right)_n}
\sqrt{\frac{2(2\lambda +1)_n(n+\lambda )}{
n!(n+2\lambda )}}.
$$ 

In order to compute the parameters $\{\gamma ^{\lambda }_{k,j}\}$
of the weight function $w$ we use Theorem ~
\ref{nouaparam}. Therefore we need to calculate the values
$\varphi ^{\lambda }_n(0,l)$ and $k^{\lambda ,l}_n$, $n\geq 1$, $l\geq 0$
where $k^{\lambda ,l}_n$ denotes the leading coefficient of 
$\varphi ^{\lambda }_n(0,l)$. The main point for these calculations 
is to notice 
that $\{\varphi ^{\lambda }_n(x,l)\}_{n\geq 0}$ is the family of
orthonormal polynomials associated to the weight function
$x^{2l}w(x)$. These polynomials are also classical objects and they
can be found for instance in \cite{DX} under the name of modified
classical polynomials. A calculation of the modified Gegenbauer
polynomials can be obtained in terms of Jacobi polynomials. These 
are orthogonal polynomials associated to parameters $\alpha ,\beta >1$
and weight function 
$$2^{-\alpha -\beta -1}B(\alpha +1,\beta +1)^{-1}
(1-x)^{\alpha }(1+x)^{\beta }$$
on $(-1,1)$ by the formula
$$P^{(\alpha ,\beta )}_n(x)=
\frac{(-1)^n}{2^nn!}
(1-x)^{-\alpha }(1+x)^{-\beta }
\frac{d^n}{dx^n}(1-x)^{\alpha +n}(1+x)^{\beta +n}.$$ 
According to \cite{DX}, Sect. 1.5.2, we have
$$\varphi ^{\lambda }_{2n}(x,l)=c_{2n}
P^{\lambda -\frac{1}{2},l-\frac{1}{2}}_n(2x^2-1)$$
and 
$$\varphi ^{\lambda }_{2n+1}(x,l)=c_{2n+1}x
P^{\lambda -\frac{1}{2},l+\frac{1}{2}}_n(2x^2-1),$$
where $c_n$ is a constant that remains to be determined. But first we can 
already notice that the above formulae give 
$\varphi ^{\lambda }_{2n+1}(0,l)=0$, so that 
$\gamma ^{\lambda }_{l,2n+1+l}=0$.
\begin{theorem}\label{gegenbauer}
For $n,l\geq 1$, 
$$\varphi ^{\lambda }_{2n}(0,l)=
(-1)^{n+1}
\sqrt{
\frac{(\lambda +1)_l}{\left(\frac{1}{2}\right)_lh^{\lambda ,l}_{2n}}}
\times \prod _{k=1}^n\frac{\lambda +l+k-1}{k}$$
and 
$$k^{\lambda ,l}_{2n}=\frac{(\lambda +l)_{2n}}{\left(l+\frac{1}{2}\right)_nn!}
\sqrt{\frac{(\lambda +1)_l}{\left(\frac{1}{2}\right)_lh^{\lambda ,l}_{2n}}},
$$
$$k^{\lambda ,l}_{2n+1}=
\frac{(\lambda +l)_{2n+1}}{\left(l+\frac{1}{2}\right)_{n+1}n!}
\sqrt{\frac{(\lambda +1)_l}{\left(\frac{1}{2}\right)_lh^{\lambda ,l}_{2n+1}}},
$$ 
where 
$$h^{\lambda ,l}_{2n}=
\frac{\left(\lambda +\frac{1}{2}\right)_n(\lambda +l)_n(\lambda +l)}
{n!\left(l+\frac{1}{2}\right)_n(\lambda +l+2n)}$$
and 
$$h^{\lambda ,l}_{2n+1}=\frac{\left(\lambda +\frac{1}{2}\right)_n
(\lambda +l)_{n+1}(\lambda +l)}
{n!\left(l+\frac{1}{2}\right)_{n+1}(\lambda +l+2n+1)}
.
$$
\end{theorem}
\begin{proof}
It is more convenient to introduce
the polynomials
$$C^{(\lambda ,l)}_{2n}(x)=\frac{(\lambda +l)_n}{\left(l+\frac{1}{2}\right)_n}
P^{(\lambda -\frac{1}{2},l-\frac{1}{2})}_n(2x^2-1),$$
$$C^{(\lambda ,l)}_{2n+1}(x)=
\frac{(\lambda +l)_{n+1}}{(l+\frac{1}{2})_{n+1}}
xP^{(\lambda -\frac{1}{2},l+\frac{1}{2})}_n(2x^2-1),$$
and again by classical results that can be found in \cite{DX}, 
we deduce
$$\begin{array}{rcl}
1\!\!\!&=&\!\!\!\displaystyle\int _{-1}^1x^{2l}
\left(\varphi ^{\lambda }_{2n}(x,l)\right)^2w(x)dx \\
 & & \\
 &=&\!\!\!c_{2n}^2\left(\frac{\left(l+\frac{1}{2}\right)_n}{(\lambda +l)_n}
\right)^2\frac{B(l+\frac{1}{2},\lambda +\frac{1}{2})}
{B(\frac{1}{2},\lambda +\frac{1}{2})}
\displaystyle\int _{-1}^1x^{2l}\left(C^{\lambda ,l}_{2n}(x)\right)^2
\frac{1}{B(l+\frac{1}{2},\lambda +\frac{1}{2})}
(1-x^2)^{\lambda -\frac{1}{2}}dx \\
 & & \\
 &=&\!\!\!c_{2n}^2\left(\frac{(l+\frac{1}{2})_n}{(\lambda +l)_n}
\right)^2\frac{B(l+\frac{1}{2},\lambda +\frac{1}{2})}
{B(\frac{1}{2},\lambda +\frac{1}{2})}h^{\lambda ,l}_{2n},
\end{array}
$$
where 
$$h^{\lambda ,l}_{2n}=
\frac{(\lambda +\frac{1}{2})_n(\lambda +l)_n(\lambda +l)}
{n!\left(l+\frac{1}{2}\right)_n(\lambda +l+2n)}.$$
Using that $\frac{B(l+\frac{1}{2},\lambda +\frac{1}{2})}
{B(\frac{1}{2},\lambda +\frac{1}{2})}=\frac{\left(\frac{1}{2}\right)_l}
{(\lambda +1)_l}$, we deduce 
$$\varphi ^{\lambda }_{2n}(x,l)=\sqrt{
\frac{(\lambda +1)_l}{\left(\frac{1}{2}\right)_lh^{\lambda ,l}_{2n}}}
\times C^{\lambda ,l}_{2n}(x).$$
The calculation of $\varphi ^{\lambda }_{2n}(0,l)$
reduces to the calculation of $C^{\lambda ,l}_{2n}(0)$
which can be easily done due to the three term relation
$$C^{\lambda ,l}_{2n+2}(x)=
\frac{\lambda +l+2n+1}{n+1}xC^{\lambda ,l}_{2n+1}(x)-
\frac{\lambda +l+n}{n+1}C^{\lambda ,l}_{2n}(x).$$
Thus we deduce 
$$C^{\lambda ,l}_{2n+2}(0)=-\frac{\lambda +l+n}{n+1}C^{\lambda ,l}_{2n}(0),
$$
and by iterating this relation and using that $C^{\lambda ,l}_{0}(0)=1$, 
we get
$$\varphi ^{\lambda }_{2n}(0,l)
=(-1)^{n+1}
\sqrt{
\frac{(\lambda +1)_l}{\left(\frac{1}{2}\right)_lh^{\lambda ,l}_{2n}}}
\times \prod _{k=1}^n\frac{\lambda +l+k-1}{k}.$$
The leading coefficient of $\varphi ^{\lambda }_{2n}(x,l)$ can be obtained 
from the corresponding formula in \cite{DX}. Thus, 
$$k^{\lambda ,l}_{2n}=
\frac{(\lambda +l)_{2n}}{\left(l+\frac{1}{2}\right)_nn!}
\sqrt{\frac{(\lambda +1)_l}{\left(\frac{1}{2}\right)_lh^{\lambda ,l}_{2n}}}
$$
and 
$$k^{\lambda ,l}_{2n+1}=
\frac{(\lambda +l)_{2n+1}}{\left(l+\frac{1}{2}\right)_{n+1}n!}
\sqrt{\frac{(\lambda +1)_l}{\left(\frac{1}{2}\right)_lh^{\lambda ,l}_{2n+1}}},
$$ 
where 
$$h^{\lambda ,l}_{2n+1}=\frac{\left(\lambda +\frac{1}{2}\right)_n
(\lambda +l)_{n+1}(\lambda +l)}
{n!\left(l+\frac{1}{2}\right)_{n+1}(\lambda +l+2n+1)}.$$
\end{proof}
Now the parameters $\{\gamma ^{\lambda }_{k,j}\}$
can be easily 
calculated by using Theorem ~\ref{gegenbauer}. Of course, 
the explicit formulae 
look too complicated to be recorded here.

\bigskip
\noindent
{\it 3.5. Asymptotic properties.}
In the classical setting of orthogonal polynomials on the unit circle there
are several remarkable asymptotic results given by Szeg\"o.
Let $\mu $ be a measure in the Szeg\"o class, 
and let $\{\varphi _n\}_{n\geq 0}$ be the family of orthonormal 
polynomials associated to $\mu $. Then, the orthonormal polynomials
have the following asymptotic properties:
\begin{equation}\label{unu}
\varphi _n\rightarrow 0
\end{equation}
and 
\begin{equation}\label{doi}
\frac{1}{\varphi _n^{\sharp }}\rightarrow \Theta _{\mu },
\end{equation}
where $\Theta _{\mu}$ is the spectral factor of $\mu $ and the convergence is 
uniform on the compact subsets of the unit disk $\DD $. The second limit 
\eqref{doi} is related to the so-called Szeg\"o limit theorems concerning
the asymptotic behaviour of Toeplitz determinants. Thus, 
$$\frac{\det T_n}{\det T_{n-1}}=\frac{1}{|\varphi _n^{\sharp }(0)|^2},
$$
where $T_n=[s_{i-j}]_{i,j=0}^n$ and $\{s_k\}_{k\in \ZZ }$ is the set
of the Fourier coefficients of $\mu $. As a consequence of the previous 
relation and \eqref{doi} we deduce Szeg\"o's first limit theorem,
\begin{equation}\label{trei}
\lim _{n\rightarrow \infty }
\frac{\det T_n}{\det T_{n-1}}=|\Theta _{\mu }(0)|^2=
\exp (\frac{1}{2\pi }\int _0^{2\pi }\log \mu '(t)dt).
\end{equation}
The second (strong) Szeg\"o limit theorem improves \eqref{trei}
by showing that 
\begin{equation}\label{patru}
\lim _{n\rightarrow \infty }
\frac{\det T_n}{g^{n+1}(\mu )}=
\exp \left(\frac{1}{\pi }\int\int _{|z|\leq 1}
|\Theta '_{\mu }(z)/\Theta _{\mu }(z)|^2d\sigma (z)\right),
\end{equation}
where $g(\mu )$ is the limit in \eqref{trei} and $\sigma $ is the 
planar Lebesgue measure. These two limits \eqref{trei}
and \eqref{patru} have an useful interpretation in terms of asymptotics 
of angles in the geometry of a stochastic process associated to $\mu $
(see \cite{GS}) and many important applications. 
We show how these results can be extended to orthogonal polynomials
on $\cP _1$. The formulae
\eqref{primarelatie} and \eqref{adouarelatie} suggest that it is 
more convenient to work in a larger algebra. This is related to the so-called
Toeplitz embedding, see \cite{DGK}, \cite{FFGK}. 
Thus, we consider the 
set $\cL $ of lower triangular arrays 
$a=[a_{k,j}]_{k,j\geq 0}$ with complex entries.
No boundedness assumption is made on these arrays. The addition in 
$\cL $ is defined by entry-wise addition and the multiplication is 
the matrix multiplication: for $a=[a_{k,j}]_{k\geq j}$, 
$b=[b_{k,j}]_{k,j\geq 0}$ two elements of $\cL $,
$$(ab)_{k,j}=\sum _{l\geq 0}a_{k,l}b_{l,j},$$
which is well-defined since the sum is finite. Thus,
$\cL $ becomes an associative, unital algebra. 

Next we associate the element $\Phi _n$ of $\cL $  
to the polynomials 
$\varphi _n(X_1,l)=\sum _{k=0}^na_{n,k}^lX_1^k$, $n,l\geq 0$, by the formula
\begin{equation}\label{fiunu}
(\Phi _n)_{k,j}=\left\{\begin{array}{lcl}
a_{n,k-j}^j & \quad & k\geq j \\
0    & \quad & k<j;
\end{array}\right.
\end{equation}
similarly, the element $\Phi ^{\sharp }_n$ of $\cL $ is associated 
to the family of polynomials 
$\varphi ^{\sharp }_n(X_1,l)=\sum _{k=0}^nb_{n,k}^lX_1^k$, $n,l\geq 0$, 
by the formula
\begin{equation}\label{fidoi}
(\Phi ^{\sharp }_n)_{k,j}=\left\{\begin{array}{lcl}
b_{n,k-j}^j &\quad  & k\geq j \\
0 & \quad  & k<j.
\end{array}\right.
\end{equation}

We notice that the spectral factor 
$\Theta _{\phi }$ of $K_{\phi }$ is an element of
$\cL $ and we assume that $\Theta _{\phi }$ belongs to the Szeg\"o class.
This implies that 
$\Phi ^{\sharp }_n$ is invertible in $\cL$ for all $n\geq 0$.
Finally, we say that a sequence $\{a_n\}\subset  \cL $
converges to $a\in \cL $ if $\{(a_n)_{k,j}\}$ 
converges to $a_{k,j}$ for all $k,j\geq 0$ (and we write $a_n\rightarrow a$).
\begin{theorem}\label{convergenta}
Let $\phi $ belong to the Szeg\"o class. Then 

\begin{equation}\label{doiunu}
\Phi _n\rightarrow 0
\end{equation}
and 
\begin{equation}\label{doidoi}
(\Phi _n^{\sharp })^{-1}\rightarrow \Theta _{\phi }.
\end{equation}
\end{theorem}

We now briefly discuss the geometric setting for the kernel 
$K_{\phi }$. By a classical result of Kolmogorov (see \cite{Pa}), 
$K_{\phi }$ is the covariance 
kernel of a stochastic process $\{f_n\}_{n\geq 0}\subset L^2(\mu )$
for some probability space $(X,\cM ,\mu )$. That is, 
$$K_{\phi }(m,n)=\int _Xf_n\overline{f}_md\mu .$$
We can suppose, without loss of generality, that $\{f_n\}_{n\geq 0}$
is total in $L^2(\mu )$ and for $p\leq q$ we introduce 
the subspaces $\cE _{p,q}$ 
given by the closure in $L^2(\mu )$ of the linear span of 
$\{f_k\}_{k=p}^q$. 
The operator angle between two spaces $\cE _1$ and $\cE _2$ of 
$L^2(\mu )$ is defined by 
$$B(\cE _1, \cE _2)=P_{\cE _1}P_{\cE _2}P_{\cE _1},$$
where $P_{\cE _1}$ is the orthogonal projection of 
$L^2(\mu )$ onto $\cE _1$. Also define 
$$\Delta (\cE _1, \cE _2)=I-B(\cE _1, \cE _2).$$

We associate to the process $\{f_n\}_{n\geq 0}$
a family of subspaces $\cH _{r,q}$ of $L^2(\mu )$
such that $\cH _{r,q}$ is the closure of the linear space generated
by $f_k$, $r\leq k\leq q$ and 
we consider 
a scale of limits:
\begin{equation}\label{alfa}
s-\lim _{r\rightarrow \infty }\Delta (\cH _{0,n},\cH _{n+1,r})=
\Delta (\cH _{0,n},\cH _{n+1,\infty })
\end{equation}
for $n\geq 0$, and then we let $n\rightarrow \infty $ and deduce
\begin{equation}\label{beta}
s-\lim _{n\rightarrow \infty }\Delta (\cH _{0,n},\cH _{n+1,\infty })=
\Delta (\cH _{0,\infty },\cap _{n\geq 0}\cH _{n,\infty }),
\end{equation} 
where $s-\lim $ denotes the strong operatorial limit. 

We then deduce analogues of the Szeg\"o limit theorems \eqref{trei}
and \eqref{patru} by expressing 
the above limits of angles in terms of determinants.
This is possible due to \eqref{deter}.
\begin{theorem}\label{szego}
Let $\phi $ belong to the Szeg\"o class.
Then
\begin{equation}\label{3?}
\frac{D_{r,q}}{D_{r+1,q}}=s_{r,r}\det \Delta (\cH _{r,r},\cH _{r+1,q})=
\frac{1}{|\varphi ^{\sharp }_{q-r}(0,r)|^2}
\end{equation}
and 
\begin{equation}\label{3??}
\lim _{q\rightarrow \infty }\frac{D_{r,q}}{D_{r+1,q}}=
s_{r,r}\det \Delta (\cH _{r,r},\cH _{r+1,\infty })=
|\Theta _{\phi }(r,r)|^2=s_{r,r}\prod _{j\geq 1}d^2_{r,r+j}.
\end{equation}
If we denote the above limit by $g_r$ and 
$$L=\lim _{n\rightarrow \infty }\prod _{0\leq k<n<j}
d_{k,j}^2>0,$$
then 
\begin{equation}\label{3???}
\lim _{n\rightarrow \infty }
\frac{D_{0,n}}{\prod _{l=0}^ng_l}
=\frac{1}{\det \Delta (\cH _{0,\infty }, \cap _{n\geq 0}\cH _{n,\infty })}
=\frac{1}{L}.
\end{equation}
\end{theorem}
\noindent
Details of the proofs can be found in \cite{BC}.

\bigskip
\noindent
{\it 3.6. Szeg\"o kernels.}
The classical theory of orthogonal polynomials is intimately related to 
some classes of analytic functions. Much of this interplay is realized by
Szeg\"o kernels. Here we expand this idea by 
providing a Szeg\"o type kernel for 
(a slight modification of) the space $\cH _0^2(\cF )$ which is viewed as an 
analogue of the Hardy class $H^2$ on the unit disk. We mention that 
another version of this idea was developed in \cite{ADD} (and recently
applied to a setting of stochastic processes indexed by vertices of 
homogeneous trees in \cite{AV}); see also \cite{Ar0},
\cite{BaG}. The difference is that $\cH _0^2(\cF )$
is larger than the space $\cU _2$ of \cite{ADD} and also, the Szeg\"o
kernel that we consider is positive definite.

We return now to the space $\cH _0^2(\cF )$ introduced in Subsection ~3.2.
Its definition involves a family of Hilbertian conditions therefore its
natural structure should be that of a Hilbert module
(we use the terminology of \cite{La}). 
Assume $\cF _n=\CC $ for all $n\geq 0$ and consider the 
$C^*$-algebra $\cD $ of bounded diagonal operators 
on $\oplus _{n\geq 0}\cF _n=l^2(\NN _0)$.
For a sequence $\{d_n\}_{n\geq 0}$ we use the notation
$$\mbox{diag}\left(\{d_n\}_{n\geq 0}\right)=
\left[
\begin{array}{ccccc}
d_0 & 0 & \hdots & & \\
0 & d_1 & \hdots & \mbox{\large{0}} & \\
\vdots  & \vdots & d_2 & & \\
& \mbox{\large{0}} & &  \ddots & \\
 & & & & 
\end{array}
\right],
$$
so that $\mbox{diag}\left(\{d_n\}_{n\geq 0}\right)$ 
belongs to $\cD $ if and only if
$\sup _{n\geq 0}|d_n|<\infty $. We now define the vector space
$$\cH ^2(\cF )=\{\Theta \in \cH _0^2(\cF )\mid 
\mbox{diag}\left(\{c_n(\Theta )^*c_n(\Theta )\}_{n\geq 0}\right)
\in \cD \}$$
and notice that $\cD $ acts linearly on $\cH ^2(\cF )$ by $\Theta D=
\left[\Theta _{k,j}d_j\right]_{k,j\geq 0}$, therefore
$\cH ^2(\cF )$ is a right $\cD$-module. Also, 
if $\Theta $, $\Psi $ belong to $\cH ^2(\cF )$, then 
$\mbox{diag}\left(\{c_n(\Psi )\}^*c_n(\Theta )\}_{n\geq 0}\right)$
belongs to $\cD $, which allows us to define
$$\langle \Theta ,\Psi \rangle =
\mbox{diag}\left(\{c_n(\Psi )^*c_n(\Theta )\}_{n\geq 0}\right),
$$
turning  $\cH ^2(\cF )$ into a Hilbert $\cD $-module.
As a Banach space with norm 
$\|\Theta \|=\|\langle \Theta ,\Theta \rangle \|^{1/2}$, 
the space $\cH ^2(\cF )$
coincides with $l^{\infty }(\NN _0,l^2(\NN _0))$, the 
Banach space of bounded sequences of elements in $l^2(\NN _0)$.
A similar construction is used in \cite{BG} for a setting of
orthogonalization with invertible squares. 

Next we introduce a Szeg\"o kernel for $\cH ^2(\cF )$.
Consider the set
$$B_1=\{\{z_n\}_{n\geq 0}\subset \CC \mid \sup _{n\geq 0}|z_n|<1\}$$
and for $z=\{z_n\}_{n\geq 0}\in B_1$, $\Theta \in  \cH ^2(\cF )$,
notice that 
$$\Theta (z)=\mbox{diag}\left(\{
\Theta _{n,n}+\sum _{k>n}\Theta _{k,n}z_{k-1}\ldots z_n\}_{n\geq 0}\right)
$$
is a well-defined element of $\cD$. Also, for $z\in B_1$, we define
$$S_z=
\left[\begin{array}{cccc}
1 & 0 & 0 & \hdots  \\
\overline{z}_0 & 1 & 0 &  \hdots  \\
\overline{z}_0\overline{z}_1 & \overline{z}_1 & 1 &  \hdots  \\
\overline{z}_0\overline{z}_1 \overline{z}_2 &
\overline{z}_1\overline{z}_2 &
\overline{z}_2 &  \hdots \\ 
\vdots & \vdots & \vdots & \ddots  
\end{array}
\right],
$$
which is an element of $\cH ^2(\cF )$
and the Szeg\"o kernel is defined on $B_1$ by the formula:
\begin{equation}\label{kszego}
S(z,w)=\langle S_w,S_z\rangle ,\quad z,w\in B_1.
\end{equation}

\begin{theorem}\label{kernel}
$S$ is a positive definite kernel on $B_1$ with the 
properties:

\smallskip
\quad $(1)$ \quad $\Theta (z)=\langle \theta ,S_z\rangle ,
\quad \Theta \in  \cH ^2(\cF ), z\in B_1$.

\smallskip
\quad $(2)$ \quad
The set $\{S_zD \mid z\in B_1, D\in \cD \}$ is total in $\cH ^2(\cF )$.
\end{theorem}
\begin{proof}
Take $z_1,\ldots z_m\in B_1$ and after reshuffling the matrix 
$\left[S(z_j,z_l)\right]_{j,l=1}^m$ can be written in the form
$$\oplus _{n\geq 0}\left[c_n(S_{z_j})^*c_n(S_{z_l})\right]_{j,l=1}^m.
$$
Since each matrix $\left[c_n(S_{z_j})^*c_n(S_{z_l})\right]_{j,l=1}^m=
\left[\begin{array}{c}
c_n(S_{z_1}) \\
\vdots \\
c_n(S_{z_m})
\end{array}
\right]^*
\left[\begin{array}{c}
c_n(S_{z_1}) \\
\vdots \\
c_n(S_{z_m})
\end{array}
\right]$
is positive, we conclude that $S$ is a positive definite kernel
on $B_1$. Property $(1)$ of $S$ follows directly from definitions.
For $(2)$ we use approximation theory in $L^{\infty }$
spaces in order to reduce the proof to the 
following statement: for every
element $\Theta \in \cH ^2(\cF )$ for which  
there exists $A\subset \NN $ such that 
$c_n(\Theta )=h\in l^2(\NN _0)$
for $n\in A$ and 
$c_n(\Theta )=0$ for $n\notin A$, and for every $\epsilon >0$, 
there exists a linear combination $L$ of elements $S_zD$, 
$z\in B_1$, $D\in \cD $, such that 
$\sup _{n\geq 0}\|c_n(\Theta -L)\|<\epsilon .$ 
This can be achieved as follows. Since the set 
$\{\phi _{w}(z)=\frac{1}{1-z\overline{w}}\mid w\in B_1\}$ is total 
in the Hardy space $H^2$ on the unit disk, we deduce that there
exist complex numbers $c_1$, $\ldots $, $c_m$ and
$w_1$, $\ldots $, $w_m$, $|w_k|<1$ for all $k=1,\ldots m$, 
such that 
$$\|h-\sum _{k=1}^mc_k\left[\begin{array}{c}
1 \\
\overline{w}_k \\
\overline{w}^2_k \\
\vdots 
\end{array}
\right]
\|<\epsilon.
$$
Then define $z_k=\{w_{k,n}\}_{n\geq 0}$ for $k=1,\ldots ,m$, 
where $w_{k,n}=w_k$ for all $n\geq 0$. So $z_k\in B_1$. Also define
$d^A_n=1$ for $n\in A$ and $d^A_n=0$ for $n\notin A$, and consider
$$L=\sum _{k=1}^mc_kS_{z_k}\mbox{diag}\left(\{d^A_n\}\right).$$
We deduce that $L\in \cH ^2(\cF) $, $c_n(L)=0$ for 
$n\notin A$, and $c_n(L)=\sum _{k=1}^mc_k\left[\begin{array}{c}
1 \\
\overline{w}_k \\
\overline{w}^2_k \\
\vdots 
\end{array}
\right]
$ for $n\in A$, 
so that 
$$\begin{array}{rcl}
\sup _{n\geq 0}\|c_n(\Theta -L)\}
& = & \sup _{n\geq 0}\|c_n(\Theta )-c_n(L)\} \\
& & \\
& = & \max \{\sup _{n\in A}\|h-c_n(L)\|,\sup _{n\notin A}\|c_n(L)\|\} \\
& & \\
& = & \|h-\sum _{k=1}^mc_k\left[\begin{array}{c}
1 \\
\overline{w}_k \\
\overline{w}^2_k \\
\vdots 
\end{array}
\right]
\|<\epsilon.
\end{array}
$$
\end{proof}

\section{Several isometric variables}
In this section we discuss orthogonal polynomials in several variables 
satisfying the isometric relations $X^+_kX_k=1$, $k=1,\ldots ,N$.
We set $\cA =\{1-X^+_kX_k\mid k=1,\ldots ,N\}$ and notice that the index set
of $\cA $ is $\FF ^+_N$. 
Also if $\phi $ is a linear functional on $\cR (\cA )$ 
then its kernel of moments is invariant under the action of $\FF ^+_N$
on itself by juxtaposition, that is,
\begin{equation}\label{sta1}
K_{\phi }(\tau \sigma ,\tau \sigma ')=K_{\phi }(\sigma ,\sigma '), \quad 
\tau ,\sigma ,\sigma '\in \FF _N^+.
\end{equation} 
In fact, a kernel $K$ satisfies \eqref{sta1} if and only if $K=K_{\phi }$ for
some linear functional on $\cR (\cA )$. Positive definite kernels
satisfying \eqref{sta1} have been already studied, see for instance 
\cite{BBCGNW} and references therein. In particular, the class
of isotropic processes on homogeneous trees give rise to positive definite
kernels for which a theory of orthogonal polynomials (Levinson recursions)
was developed in \cite{BBW}. Here we discuss in more details another class of
kernels satisfying \eqref{sta1} which was considered, for instance, 
in \cite{Fr}.

\bigskip
\noindent
{\it 4.1. Cuntz-Toeplitz relations.}
Consider the class of positive definite kernels with property 
\eqref{sta1} and such that
\begin{equation}\label{sta2}
K(\sigma ,\tau )=0 \quad \mbox{ if there is no $\alpha \in \FF _N^+$
such that $\sigma =\tau \alpha $ or $\tau =\sigma \alpha $}.
\end{equation}

We showed in \cite{CJ}
that $K$ has properties \eqref{sta1} and \eqref{sta2}
if and only if $K=K_{\phi }$ for some q-positive functional 
on $\cR (\cA _{CT})$, where 
$\cA _{CT}=\{1-X^+_kX_k\mid k=1,\ldots ,N\}\cup
\{X^+_kX_l, k,l=1,\ldots ,N, k\ne l\}$. The relations in 
$\cA _{CT}$ are defining the Cuntz-Toeplitz algebra (see
\cite{EK} for details). The property \eqref{sta2} shows that $K$ is 
quite sparse, therefore it is expected to be easy to analyse such a kernel.
Still, there are some interesting aspects related to this class of kernels, 
some of which we discuss here.

Let $\phi $ be a strictly q-positive kernel on $\cR (\cA _{CT})$ and let
$K_{\phi }$ be the associated kernel of moments. Since
the index set of $\cA _{CT}$ is still $\FF ^+_N$, and this is totally 
ordered by the lexicographic order, we can use the results described in 
Subsection ~3.1 and associate to $K_{\phi }$ a family 
$\{\gamma _{\sigma , \tau }\}_{\sigma \prec \tau }$
of complex numbers with $|\gamma _{\sigma , \tau }|<1$,
uniquely determining $K_{\phi }$ by relations of type 
\eqref{gdil}. It was noticed in \cite{CJ} that $K_{\phi }$ has
properties \eqref{sta1} and \eqref{sta2} if and only if 
$\gamma _{\tau \sigma ,\tau \sigma '}=\gamma _{\sigma ,\sigma '}$ 
and
$
\gamma _{\sigma ,\tau }=0$ if there is no $\alpha \in \FF _N^+$
such that $\sigma =\tau \alpha $ or $\tau =\sigma \alpha $.
The main consequence of these relations is that $K_{\phi }$ is uniquely
determined by $\gamma _{\sigma }=\gamma _{\emptyset ,\sigma }$, 
$\sigma \in \FF _N^+-\{\emptyset \}$. We define $d_{\sigma }=
(1-|\gamma _{\sigma }|^2)^{1/2}$.
The orthogonal polynomials associated to $\phi $ satisfy the following
recurrence relations (see \cite{CJ} for details): 
$\varphi _{\emptyset }=\varphi ^{\sharp }_{\emptyset }=
s^{-1/2}_{\emptyset ,\emptyset }$ and for $k\in \{1,\ldots ,N\}$,
$\sigma \in \FF _N^+$, 
\begin{equation}\label{bszego}
\varphi _{k\sigma }=\frac{1}{d_{k\sigma }}
(X_k\varphi _{\sigma }-\gamma _{k\sigma } 
\varphi ^{\sharp}_{k\sigma -1}),
\end{equation}
\begin{equation}\label{bsarp}
\varphi ^{\sharp}_{k\sigma }=\frac{1}{d_{k\sigma }}
(-\overline{\gamma }_{k\sigma }X_k\varphi _{\sigma }+ 
\varphi ^{\sharp}_{k\sigma -1}).
\end{equation}
The results corresponding to Theorem ~\ref{convergenta} and 
Theorem ~\ref{szego} can be easily obtained (see \cite{BC}), 
but the constructions around the Szeg\"o kernel are more interesting
in this situation. Thus, there is only one Hilbertian condition involved
in the definition of $\cH _0(\cF )$ in this case. In fact, it is easy
to see that $\cH _0(\cF )$ can be identified with the full Fock space 
$l^2(\FF ^+_N)$, the $l^2$ space over $\FF ^+_N$. Now, concerning 
evaluation of elements of 
$\cH _0(\cF )$, if we are going to be consistent with the point 
of view that the ``points for evaluation'' come from
the unital homomorphisms of the polynomial algebra   
inside $\cH _0(\cF )$, then we have to consider an infinite 
dimensional Hilbert space $\cE $ and the set
$$B_1(\cE )=\{Z=(Z_1,\ldots ,Z_N)\in \cL (\cE )^N 
\mid \sum _{k=1}^NZ_kZ^*_k<I\}.
$$
For $\sigma =i_1\ldots i_k\in \FF ^+_N$
we write $Z_{\sigma }$ instead of $Z_{i_1}\ldots Z_{i_k}$. 
Then we define for $\Theta \in l^2(\FF ^+_N)\otimes \cE $ and 
$Z\in B_1(\cE )$, 
$$\Theta (Z)=\sum _{\sigma \in \FF ^+_N}Z_{\sigma }\Theta _{\sigma },$$
which is an element of the set $\cL (\cE )$ of bounded linear operators
on the Hilbert space $\cE $. 
Next,
for $Z\in B_1(\cE )$ we define $S_Z:\cE \rightarrow l^2(\FF ^+_N)\otimes \cE $
by the formula:
$$S_Zf=\sum _{\sigma \in \FF ^+_N}e_{\sigma }\otimes (Z_{\sigma })^*f,
\quad f\in \cE .$$
Then $S_Z\in \cL (\cE ,l^2(\FF ^+_N)\otimes \cE )$ and we 
can finally introduce the Szeg\"o kernel on $B_1(\cE )$ by the formula:
$$S(Z,W)=S_Z^*S_W,\quad Z,W\in B_1(\cE ).$$

\begin{theorem}\label{doiker}
$S$ is a positive definite kernel on $B_1(\cE )$ with the 
properties:

\smallskip
\quad $(1)$ \quad $\Theta (z)=S_Z^*\Theta ,
\quad \Theta \in l^2(\FF ^+_N)\otimes \cE , Z\in B_1(\cE )$.

\smallskip
\quad $(2)$ \quad
The set $\{S_Zf \mid Z\in B_1(\cE ), f\in \cE \}$ is total in 
$l^2(\FF ^+_N)\otimes \cE $.
\end{theorem}
\begin{proof}
The fact that $S$ is positive definite and $(1)$ are immediate.
More interesting is $(2)$ and we reproduce here the proof 
given in \cite{CJ}.
Let $f=\{f_{\sigma }\}_{\sigma \in \FF ^+_N}$
be an element of $l^2(\FF ^+_N)\otimes \cE $ orthogonal
to the linear space generated by $\{S_Zf \mid Z\in B_1(\cE ), f\in \cE \}$.
Taking $Z=0$, we deduce that $f_{\emptyset }=0$. Next, we claim
that for each $\sigma \in \FF ^+_N-\{\emptyset \}$ there exist
$$Z_l=(Z_1^l,\ldots ,Z_N^l)\in \cB _1(\cE ), \quad 
l=1,\ldots ,2|\sigma |,$$
such that 
$$\mbox{range}\left[
\begin{array}{ccc}
Z_{\sigma }^{*1} & \ldots & Z_{\sigma }^{*2|\sigma |}
\end{array}\right]=\cE,$$
and 
$$Z_{\tau }^l=0\quad \mbox{for all}\quad 
\tau \ne \sigma , \quad |\tau |\geq |\sigma |, \quad l=1,\ldots ,2|\sigma |.$$

Once this claim is proved, a simple inductive argument gives
$f=0$.
In order to prove the claim we need the following construction.
Let $\{e^n_{ij}\}_{i,j=1}^n$ be the matrix units of the algebra
$M_n$ of $n\times n$ matrices. Each 
$e^n_{ij}$ is an $n\times n$ matrix consisting of $1$
in the $(i,j)th$ entry and zeros elsewhere.
For a Hilbert space $\cE _1$ we define $E^n_{ij}=e^n_{ij}\otimes 
I_{\cE _1}$ and we notice that 
$E^n_{ij}E^n_{kl}=\delta _{jk}E^n_{il}$ and $E^{*n}_{ji}=E^n_{ij}$.
Let $\sigma =i_1\ldots i_k$ so that $\cE =\cE _1^{\oplus 2|\sigma |}$
for some Hilbert space $\cE _1$ (here we use in an essential way 
the assumption 
that $\cE $ is of infinite dimension).
Also, for $s=1,\ldots ,N$, we define 
$J_s=\{l\in \{1,\ldots ,k\}\mid i_{k+1-l}=s\}$
and 
$$
Z^{*p}_{s }=\frac{1}{\sqrt{2}}
\sum _{r\in J_s}E^{2|\sigma |}_{r+p-1,r+p},\quad 
s=1,\ldots ,N, \quad p=1,\ldots ,|\sigma |.$$
We can show that for each $p\in \{1,\ldots ,|\sigma |\}$,
\begin{equation}\label{kiwi}
Z^{*p}_{\sigma }=\frac{1}{{\sqrt{2^k}}}E^{2|\sigma |}_{p,k+p},
\end{equation}
\begin{equation}\label{mango}
Z_{\tau }^p=0\quad \mbox{for} \quad 
\tau \ne \sigma ,\quad |\tau |\geq |\sigma |.
\end{equation}
We deduce 
$$\begin{array}{rcl}
\sum _{s=1}^NZ_s^pZ_s^{*p}&=&
\frac{1}{2}\sum _{s=1}^N\sum _{r\in J_s}E^{2|\sigma |}_{r+p,r+p-1}
E^{2|\sigma |}_{r+p-1,r+p} \\
 & & \\
 &=&\frac{1}{2}\sum _{s=1}^N\sum _{r\in J_s}E^{2|\sigma |}_{r+p,r+p} \\
 & & \\
 &=&\frac{1}{2}\sum _{r=1}^kE^{2|\sigma |}_{r+p,r+p}<I,
\end{array}
$$
hence $Z^p\in B_1(\cE )$ for each $p=1,\ldots ,|\sigma |$.
For each word $\tau =j_1\ldots j_k\in \FF _N^+-\{\emptyset \}$
we deduce by induction that 
\begin{equation}\label{apple}
Z^{*p}_{j_k}\ldots 
Z^{*p}_{j_1}=\frac{1}{\sqrt{2^k}}
\sum _{r\in A_{\tau }}E^{2|\sigma |}_{r+p-1,r+p+k-1},
\end{equation}
where 
$A_{\tau }=\cap _{p=0}^{k-1}(J_{j_{k-p}}-p)\subset \{1,\ldots ,N\}$
and $J_{j_{k-p}}-p=\{l-p\mid l\in J_{i_{k-p}}\}$.

We show that $A_{\sigma }=\{1\}$ and $A_{\tau }=\emptyset $
for $\tau \ne \sigma $. Let $q\in A_{\tau }$. Therefore, for any 
$p\in \{0,\ldots ,k-1\}$ we must have $q+p\in J_{j_{k-p}}$
or $i_{k+1-q-p}=j_{k-p}$. For $p=k-1$ we deduce 
$j_1=i_{2-q}$ and since $2-q\geq 1$, it follows that $q\leq 1$. 
Also $q\geq 1$, therefore the only element that 
can be in $A_{\tau }$ is $q=1$, in which case we must have
$\tau =\sigma $. Since $l\in J_{i_{k+1-l}}$
for each $l=1, \ldots ,k-1$, hence $A_{\sigma }=\{1\}$
and $A_{\tau }=\emptyset $ for $\tau \ne \sigma $.
Formula \eqref{apple} implies \eqref{kiwi}. In a similar manner
we can construct a family $Z^p$, $p=|\sigma |+1, \ldots ,2|\sigma |,$
such that 
$$Z^{*p}_{\sigma }=\frac{1}{\sqrt{2^k}}E^{2|\sigma |}_{p+k,p},$$
and 
$$
Z_{\tau }^p=0\quad \mbox{for} \quad 
\tau \ne \sigma ,\quad |\tau |\geq |\sigma |.
$$
Thus, for $s=1,\ldots ,N$, we define 
$K_s=\{l\in \{1,\ldots ,k\}\mid i_k=s\}$
and 
$$Z^{*p}_{s}=\frac{1}{\sqrt{2}}
\sum _{r\in K_s}E^{2|\sigma |}_{r+p-k,r+p-k-1},\quad 
s=1,\ldots ,N, \quad p=|\sigma |+1,\ldots ,2|\sigma |.$$

Now, 
$$\left[
\begin{array}{ccc}
Z_{\sigma }^{*1} & \ldots & 
Z_{\sigma }^{*2|\sigma |}
\end{array}
\right]=
\frac{1}{\sqrt{2^k}}\left[
\begin{array}{cccccc}
E^{2|\sigma |}_{1,k+1} & \ldots & 
E^{2|\sigma |}_{k,2k} & E^{2|\sigma |}_{k+1,1}
& \ldots & E^{2|\sigma |}_{2k,k}
\end{array}
\right],
$$
whose range is $\cE $. 
This concludes the proof. 
\end{proof}

It is worth noticing that property $(2)$ of $S$ is no 
longer true if $\cE $ is finite dimensional. In fact, for 
$\cE $ of dimension one the set 
$\{S_Zf \mid Z\in B_1(\cE ), f\in \cE \}$ is total in the 
symmetric Fock space of $\CC ^N$ (see \cite{Ar}).

\bigskip
\noindent
{\it 4.2. Kolmogorov decompositions and Cuntz relations.}
This is a short detour from orthogonal polynomials, in order to show a 
construction of bounded operators satisfying the Cuntz-Toeplitz 
and Cuntz relations, based on parameters 
$\{\gamma _{\sigma }\}_{\sigma \in \FF ^+_N-\{\emptyset \}}$
associated to a positive definite kernel with properties \eqref{sta1}
and \eqref{sta2}.

First we deal with the Kolmogorov decomposition of a positive definite 
kernel. This is a more abstract version of the result of Kolmogorov
already alluded to in Subsection ~3.5. For a presentation of the general
result and some applications, see \cite{EK}, \cite{Pa}.
Here we consider $K:\NN _0\times \NN _0\rightarrow \CC $ a positive definite 
kernel and let $\{\gamma _{k,j}\}$ be the family of parameters associated
to $K$ as in Subsection ~3.1. In addition we assume $K(j,j)=1$ for $j\geq 0$.
This is not a real loss of generality and it simplifies some calculations.
We also assume $|\gamma _{k,j}|<1$ for all $k<j$. Then we introduce 
for $0\leq k<j$ the operator $V_{k,j}$ on $l^2(\NN _0)$
defined by the formula
$$V_{k,j}=
(J(\gamma _{k,k+1})\oplus 1_{n-1})
(1\oplus J(\gamma _{k,k+2})\oplus 1_{n-2})\ldots 
(1_{n-1}\oplus J(\gamma _{k,j}))\oplus 0$$
and we notice that 
\begin{equation}\label{41}
W_k=s-\lim _{j\rightarrow \infty }V_{k,j}
\end{equation}
is a well-defined isometric operator on $l^2(\NN _0)$ for every $k\geq 0$.
If we define
$V(0)=I/\CC $ and $V(k)=W_0W_1\ldots W_{k-1}/\CC $ for 
$k\geq 1$, then we obtain 
the following result from \cite{Co2}.

\begin{theorem}\label{kolmo}
The map $V:\NN _0\rightarrow l^2(\NN _0)$ is the Kolmogorov
decomposition of the kernel $K$< in the sense that

\smallskip
\quad $(1)$ \quad $K(j,l)=\langle V(l),V(j)\rangle ,
\quad j,l\in \NN _0$.

\smallskip
\quad $(2)$ \quad
The set $\{V(k) \mid k\in \NN _0 \}$ is total in 
$l^2(\NN _0)$.
\end{theorem}
It is worth noticing that we can write explicitly the matrix of $W_k$:
$$\left[\begin{array}{cccc}
\gamma _{k,k+1} & d_{k,k+1}\gamma _{k,k+2}
&  
d_{k,k+1}d_{k,k+2}\gamma _{k,k+3} &\ldots \\
d_{k,k+1} & -\overline{\gamma }_{k,k+1}\gamma _{k,k+2} & 
-\overline{\gamma }_{k,k+1}d_{k,k+2}\gamma _{k,k+3} &\ldots \\ 
0 & d_{k,k+2} & -\overline{\gamma} _{k,k+2}\gamma _{k,k+3} &\ldots \\
 \vdots & 0 & d_{k,k+3} & \\ 
 & \vdots & 0 & \ddots \\
& & \vdots &  
\end{array}
\right].
$$
We also see that Theorem ~3.7 and Theorem ~4.1 produce various kinds of
Kolmogorov decompositions for the corresponding Szeg\"o kernels. Based on a 
remark in \cite{Co0}, we use Theorem \ref{kolmo} in order to obtain some 
large families of bounded operators satisfying Cuntz-Toeplitz and Cuntz 
relations. Thus, we begin with a positive definite 
kernel $K$ with properties 
\eqref{sta1}
and \eqref{sta2}.
For simplicity we also assume $K(\emptyset ,\emptyset )=1$ and let
$\{\gamma _{\sigma }\}_{\sigma \in \FF ^+_N-\{\emptyset \}}$ 
be the family of corresponding parameters. In order to be in tune with the 
setting of this paper, we assume $|\gamma _{\sigma }|<1$
for all $\sigma \in \FF ^+_N-\{\emptyset \}$. Motivated by the 
construction in Theorem \ref{kolmo}, we denote  
by $c_{\sigma }(W_0)$, 
$\sigma \in \FF _N^+-\{\emptyset \}$, the columns of the operator $W_0$.
Thus, 
$$
c_1(W_0)=\left[\begin{array}{c}
\gamma _1 \\
d_{1} \\
0 \\
\vdots 
\end{array}\right]
\quad \mbox{and}
\quad 
c_{\sigma }(W_0)=\left[\begin{array}{c}
d_{1}\ldots d_{\sigma -1}\gamma _{\sigma } \\
-\overline{\gamma }_1d_{2}\ldots 
d_{\sigma -1}\gamma _{\sigma } \\
\vdots \\
-\overline{\gamma }_{\sigma -1}\gamma _{\sigma } \\
d_{\sigma } \\
0 \\
\vdots 
\end{array}\right]
\quad \mbox{for} \quad 1\prec \sigma .
$$
We now define the isometry $U(k)$ on $l^2(\FF ^+_N)$ by the formula:
$$U(k)=\left[c_{k\tau }(W_0)\right]_{\tau \in \FF ^+_N}, \quad 
k=1,\ldots ,N.$$
\begin{theorem}\label{ctc}
$(a)$ The family $\{U_1, \ldots ,U_N\}$ satifies the Cuntz-Toeplitz relations:
$U^*_kU_l=\delta _{k,l}I$, $k,l=1,\ldots ,N$.

\noindent
$(b)$  The family $\{U_1, \ldots ,U_N\}$ satifies the Cuntz relations
$U^*_kU_l=\delta _{k,l}I$, $k,l=1,\ldots ,N$ and 
$\sum _{k=1}^NU_kU^*_k=I$, if and only if
\begin{equation}\label{cu}
\prod _{\sigma \in \FF ^+_N-\{\emptyset \}}d_{\sigma }=0.
\end{equation}
\end{theorem}
\begin{proof}
$(a)$ follows from the fact that $W_0$ is an isometry. In order to prove
$(b)$ we need a characteriztion of those $W_0$ which are unitary.
Using Proposition ~1.4.5 in \cite{Co}, we deduce that $W_0$ is unitary 
if and only if \eqref{cu} holds.
\end{proof}
It is worth mentioning that if we define 
$V(\emptyset )=I/\CC $ and $V(\sigma )=U(\sigma )/\CC $, 
$\sigma \in \FF ^+_N-\{\emptyset \}$, where 
$U(\sigma )=U(i_1)\ldots U(i_k)$ provided that $\sigma =i_1\ldots i_k$, 
then $V$ is the Kolmogorov decomposition of the kernel $K$. A second remark
here is that the condition \eqref{cu} is exactly the opposite of the 
condition for $K$ being in the Szeg\"o class. Indeed, it is easily seen that 
for a positive definite kernel with properties
\eqref{sta1}
and \eqref{sta2}, the condition \eqref{36} is equivalent with 
$\prod _{\sigma \in \FF ^+_N-\{\emptyset \}}d_{\sigma }>0$.

\section{Several hermitian variables}
The theory corresponding to this case should be an analogue of the 
theory of orthogonal polynomials on the real line. We set
$\cA =\{Y_k-Y^+_k\mid k=1,\ldots ,N\}$ and 
$\cA '=\cA \cup \{Y_kY_l-Y_lY_k
\mid k,l=1,\ldots ,N\}$, and notice that $\cR (\cA )=\cP _N$. 
Also, $\cR (\cA ')$ is isomorphic to the symmetric algebra over $\CC ^N$.
Orthogonal polynomials associated to $\cR (\cA ')$, that is, orthogonal
polynomials in several commuting variables were studied intensively in 
recent years, see \cite{DX}. In this section we analyse the noncommutative 
case. The presentation follows \cite{Co1}.

Let $\phi $ be a strictly q-positive functional on 
$\cT _N(\cA _2)$ and assume for some simplicity that 
$\phi $ is unital, $\phi (1)=1$. The index set of $\cA $
is $\FF _N^+$ and let $\{\varphi _{\sigma }\}_{\sigma \in \FF _N^+}$
be the orthonormal  polynomials associated to $\phi $.
We notice that for any $P,Q\in \cP _N$,
\begin{equation}\label{simetrie}
\begin{array}{rcl}
\langle X_kP,Q\rangle _{\phi }&=&\phi (Q^+X_kP) \\
 & & \\
&=&\phi (Q^+X^+_kP) \\
 & & \\
&=&\langle P,X_kQ\rangle _{\phi },
\end{array}
\end{equation}
which implies that the kernel of moments satisfies the relation 
$s_{\alpha \sigma ,\tau }=s_{\sigma ,I (\alpha )\tau}$ for 
$\alpha ,\sigma ,\tau \in \FF _N^+$,
where $I $ denotes the involution on  
$\FF _N^+$ given by $I(i_1\ldots i_k)=i_k\ldots i_1$.
This can be viewed as a Hankel type condition, and we already noticed that
even in the one dimensional case the parameters $\{\gamma _{k.j}\}$ of the 
kernel of moments of a Hankel type are more difficult to be used. 
Therefore, we try to deduce three-terms relations for the orthonormal 
polynomials. A matrix-vector notation already used in the 
commutative case, turns out to be quite useful. Thus, 
for $n\geq 0$, we define 
$P_n=\left[\varphi _{\sigma }\right]_{|\sigma |=n},$
$n\geq 0$, and $P_{-1}=0$.

\begin{theorem}\label{T3}
There exist matrices $A_{n,k}$ and $B_{n,k}$ such that 
\begin{equation}\label{rere}
X_kP_n=
P_{n+1}B_{n,k}+
P_nA_{n,k}+
P_{n-1}B^*_{n-1,k}, \quad k=1,\ldots ,N, n\geq 0.
\end{equation}
\end{theorem}
Each matrix $A_{n,k}$ is a selfadjoint $N^n\times N^n$ matrix, while
each $B_{n,k}$ is an $n^{n+1}\times N^n$ matrix such that
$$B_n=\left[\begin{array}{ccc}
B_{n,1} & \ldots & B_{n,N}
\end{array}
\right]
$$
is an upper triangular invertible matrix for every $n\geq 0$.
For $n=-1$, $B_{-1,k}=0$, $k=1,\ldots ,N$.
The fact that $B_n$ is upper triangular comes from the order that we use on 
$\FF ^+_N$. The invertibility of $B$ is a consequence of the fact
that $\phi $ is strictly q-positive and appears to be a basic translation of 
this information. It turns out that there are no other restrictions on the 
matrices $A_{n,k}$, $B_{n,k}$ as shown by the following Favard type result.

\begin{theorem}\label{T4}
Let $\varphi _{\sigma }=\sum _{\tau \preceq \sigma }
a_{\sigma ,\tau}X_{\tau }$, $\sigma \in \FF _N^+$,
be elements in $\cP _N$ such that  $\varphi _{\emptyset }=1$ and 
$a_{\sigma ,\sigma }>0$.
Assume that there exists a family 
$\{A_{n,k}, B_{n,k}\mid n\geq 0, \,\,k=1,\ldots ,N\}$, 
of matrices such that $A^*_{n,k}=A_{n,k}$ and  
$B_n=\left[\begin{array}{ccc}
B_{n,1} & \ldots & B_{n,N}
\end{array}
\right]$ is  an upper triangular invertible matrix for every $n\geq 0$. Also
assume that 
\begin{equation}\label{patrudoi}
X_k\left[\varphi _{\sigma }\right]_{|\sigma |=n}=
\left[\varphi _{\sigma }\right]_{|\sigma |=n+1}B_{n,k}+
\left[\varphi _{\sigma }\right]_{|\sigma |=n}A_{n,k}+
\left[\varphi _{\sigma }\right]_{|\sigma |=n-1}B^*_{n-1,k},
k=1,\ldots ,N, n\geq 0,
\end{equation}
where $\left[\varphi _{\sigma }\right]_{|\sigma |=-1}=0$ and 
$B_{-1,k}=0$ for $k=1,\ldots ,N$. 
Then there exists a unique strictly positive functional $\phi $
on $\cR (\cA )$ such that 
$\{\varphi _{\sigma }\}_{\sigma \in \FF _N^+}$
is the family of orthonormal  polynomials associated to $\phi $.
\end{theorem}

There is a family of Jacobi matrices associated to the three-term relation
in the following way.  
For $P\in \cR (\cA )(=\cP _N)$, define 
\begin{equation}\label{rep}
\Psi _{\phi }(P)\varphi _{\sigma }=P\varphi _{\sigma }.
\end{equation}
Since the kernel of moments has the Hankel type structure
mentioned above, it follows that each $\Psi _{\phi }(P)$ 
is a symmetric operator
on the Hilbert space $\cH _{\phi }$ with dense domain $\cD$, the linear space
generated by the polynomials 
$\varphi _{\sigma }$, $\sigma \in \FF ^+_N$. Moreover, 
for $P,Q \in \cP _N$, 
$$\Psi _{\phi }(PQ)=\Psi _{\phi }(P)\Psi _{\phi }(Q),$$
and $\Psi _{\phi }(P)\cD \subset \cD$,  
hence $\Psi _{\phi }$ is an unbounded representation of
$\cP _n$ (the GNS representation associated to $\phi $).
Also,  $\phi (P)=
\langle \Psi _{\phi }(P)1,1\rangle _{\phi }$ for $P\in \cP _N$.
We distinguish the operators $\Psi _k=\Psi _{\phi }(Y_k)$,
$k=1,\ldots ,N$, since
$\Psi _{\phi }(\sum _{\sigma \in \FF ^+_N}c_{\sigma }Y_{\sigma })
=\sum _{\sigma \in \FF ^+_N}c_{\sigma }\Psi _{\phi ,\sigma }$, 
where we use the notation 
$\Psi _{\phi ,\sigma }=\Psi _{i_1}\ldots \Psi _{i_k}$
for $\sigma =i_1\ldots i_k$.
Let $\{e_1,\ldots ,e_N\}$ be the standard basis
of $\CC ^N$ and define the unitary operator $W$ from $l^2(\FF ^+_N)$
onto $\cH _{\phi }$ such that 
$W(e_{\sigma })=\varphi _{\sigma }$, $\sigma \in \FF ^+_N$.
We see that $W^{-1}\cD $ is the linear space $\cD _0$
generated by $e_{\sigma }$, $\sigma \in \FF ^+_N$, 
so that we can define
$$J_k=W^{-1}\Psi _{\phi ,k}W,\quad k=1,\ldots ,N,$$ 
on $\cD _0$. Each $J_k$ is a symmetric operator on $\cD _0$
and by Theorem ~\ref{T3}, the matrix of (the closure of) $J_k$ with respect
to the orthonormal basis 
$\{e_{\sigma }\}_{\sigma \in \FF ^+_N}$ is
$$J_k=\left[
\begin{array}{cccc}
A_{0,k} & B^*_{0,k} & 0 & \ldots \\
 & & & \\
B_{0,k} & A_{1,k} & B^*_{1,k} & \\
 & & & \\ 
0 & B_{1,k} & A_{2,k} & \ddots \\
 & & & \\ 
\vdots & & \ddots & \ddots 
\end{array}
\right].$$
We call $(J_1,\ldots ,J_N)$ a Jacobi $N$-family on $\cD _0$. It is somewhat
unexpected that the usual conditions on $A_{n,k}$ and $B_{n,k}$ insure a
joint model of a Jacobi family in the following sense.
\begin{theorem}\label{Jacobi}
Let $(J_1,\ldots ,J_N)$ a Jacobi $N$-family on $\cD _0$ such that 
$A^*_{n,k}=A_{n,k}$ and 
$B_n=\left[\begin{array}{ccc}
B_{n,1} & \ldots & B_{n,N}
\end{array}\right]$ 
is an upper triangular invertible matrix for every $n\geq 0$.
Then there exists a unique strictly q-positive
functional $\phi $ on $\cP _N$ with associated orthonormal
polynomials $\{\varphi _{\sigma }\}_{\sigma \in \FF ^+_N}$ such that the map
$$W(e_{\sigma })=\varphi _{\sigma },\quad \sigma \in \FF ^+_N,$$
extends to a unitary operator from $l^2(\FF ^+_N)$ onto $\cH _{\phi }$ and 
$$J_k=W^{-1}\Psi _{\phi ,k}W, \quad k=1,\ldots ,N.$$
\end{theorem}
\begin{proof}
First the Favard type 
Theorem ~\ref{T4} gives a unique strictly q-positive functional
$\phi $  on $\cP_N$ such that 
its orthonormal polynomials satisfy the three-term 
relation associated to the given Jacobi family, and then
the GNS construction will produce the required
$W$ and $\Psi _{\phi ,k}$, as explained above.
\end{proof}
One possible application of these families of Jacobi matrices involves
some classes of random walks on $\FF ^+_N$. Figure ~4 illustrates 
an example for $N=2$ and more details are planned to be presented in 
\cite{BanC}.

\begin{figure}[h]
\setlength{\unitlength}{2700sp}%
\begingroup\makeatletter\ifx\SetFigFont\undefined%
\gdef\SetFigFont#1#2#3#4#5{%
  \reset@font\fontsize{#1}{#2pt}%
  \fontfamily{#3}\fontseries{#4}\fontshape{#5}%
  \selectfont}%
\fi\endgroup%
\begin{picture}(6466,2200)(218,-1700)
{\thinlines
\put(601,-361){\circle{150}}
}%
{\put(2401,-361){\circle{150}}
}%
{\put(301,-1261){\circle{150}}
}%
{\put(901,-1261){\circle{150}}
}%
{\put(2101,-1261){\circle{150}}
}%
{\put(2701,-1261){\circle{150}}
}%
{\put(1501,539){\circle{150}}
}%
{\put(5401,539){\circle{150}}
}%
{\put(4501,-361){\circle{150}}
}%
{\put(6301,-361){\circle{150}}
}%
{\put(4201,-1261){\circle{150}}
}%
{\put(4801,-1261){\circle{150}}
}%
{\put(6001,-1261){\circle{150}}
}%
{\put(6601,-1261){\circle{150}}
}%
{\put(1501,539){\line(-1,-1){900}}
}%
{\put(601,-361){\line(-1,-3){300}}
}%
{\put(2401,-361){\line(-5,-2){2120.690}}
}%
{\put(2401,-361){\line(-5,-3){1500}}
}%
{\put(1501,539){\line( 1,-1){900}}
}%
{\put(5401,539){\line(-1,-1){900}}
}%
{\put(5401,539){\line( 1,-1){900}}
}%
{\put(4501,-361){\line(-1,-3){300}}
}%
{\put(4501,-361){\line( 1,-3){300}}
}%
{\put(6301,-361){\line(-5,-2){2120.690}}
}%
{\put(6301,-361){\line(-5,-3){1500}}
}%
{\put(6301,-361){\line(-1,-3){300}}
}%
{\put(901,-1261){\line( 0, 1){ 75}}
}%
{\put(901,-1261){\line(-1,-1){450}}
}%
{\put(901,-1261){\line(-1,-2){225}}
}%
{\put(2101,-1261){\line(-6,-5){531.148}}
}%
{\put(2101,-1261){\line(-5,-6){375}}
}%
{\put(2101,-1261){\line(-1,-2){225}}
}%
{\put(2701,-1261){\line(-6,-5){531.148}}
}%
{\put(2701,-1261){\line(-5,-6){375}}
}%
{\put(2701,-1261){\line(-1,-2){225}}
}%
{\put(2701,-1261){\line(-1,-3){150}}
}%
{\put(4201,-1261){\line(-1,-2){225}}
}%
{\put(4201,-1261){\line( 2,-5){150}}
}%
{\put(301,-1261){\line(-1,-3){150}}
}%
{\put(4201,-1261){\line( 2,-3){300}}
}%
{\put(4801,-1261){\line(-1,-1){300}}
}%
{\put(4801,-1261){\line(-1,-2){150}}
}%
{\put(4801,-1261){\line(-1,-3){ 150}}
}%
{\put(4801,-1261){\line( 1,-4){ 75}}
}%
{\put(6001,-1261){\line(-2,-1){600}}
}%
{\put(6001,-1261){\line(-3,-2){450}}
}%
{\put(6001,-1261){\line(-3,-4){225}}
}%
{\put(6001,-1261){\line(-1,-4){ 75}}
}%
{\put(6001,-1261){\line( 1,-4){ 75}}
}%
{\put(6601,-1261){\line(-3,-1){450}}
}%
{\put(6601,-1261){\line(-2,-1){300}}
}%
{\put(6601,-1261){\line(-1,-3){150}}
}%
{\put(6601,-1261){\line(-1,-2){150}}
}%
{\put(6601,-1261){\line(-1,-4){75}}
}%
{\put(6601,-1261){\line( 2,-5){150}}
}%

\end{picture}

\caption{\mbox{ Random walks associated to a Jacobi family, $N=2$}}
\end{figure}

We conclude our discussion of orthogonal polynomials on hermitian
variables by introducing a Szeg\"o kernel that should
be related to orthogonal polynomials on $\cP _N$. Thus, we consider
the Siegel upper half-space of a Hilbert space $\cE $ by 
$$H_+(\cE )=
\{\left(W_1\ldots W_N\right)\in \cL (\cE )^N\mid W_1W^*_1+\ldots 
+W_{N-1}W^*_{N-1}<\frac{1}{2i}(W_N-W^*_N)\}.
$$
We can establish a connection between $B_1(\cE )$ and $H_+(\cE )$ similar to 
the well-known connection between the unit disk and the upper half plane 
of the complex plane. Thus, we define the Cayley transform
by the formula 
$$C(Z)=
((I+Z_N)^{-1}Z_1,\ldots ,(I+Z_N)^{-1}Z_{N-1},
i(I+Z_N)^{-1}(I-Z_N)),$$
which is well-defined for  
$Z=\left(Z_1,\ldots ,Z_N\right)\in B_1(\cE )$ since $Z_k$ must be a strict 
contraction
($\|Z_k\|<1$) for every $k=1,\ldots ,N$. 
In addition, $C$ establishes a one-to-one correspondence
from $B_1(\cE )$ onto $H_+(\cE )$. The Szeg\"o kernel on $B_1(\cE )$ can 
be transported on $H_+(\cE )$ by the Cayley transform. Thus, we introduce 
the Szeg\"o kernel on $H_+(\cE )$ by the formula:
$$S(W,W')=F^*_WF_{W'}, \quad W,W' \in H_+(\cE ),$$
where $F_W=2\mbox{diag}((-i+W^*_N))S_{C^{-1}(W)}$.
Much more remains to be done in this direction. For instance, some 
classes of orthogonal polynomials of Jacobi type and their 
generating functions are considered in \cite{BanC}.

Finally, 
we mention that there are examples of polynomial relations for which there
are no orthogonal polynomials. Thus, consider
$$\cA =\{X^+_k-X_k \mid k=1,\ldots ,N\}\cup
\{X_kX_l+X_lX_k \mid k,l=1,\ldots ,2N\},
$$
then $\cR (\cA )\simeq \Lambda (\CC ^N)$, the exterior algebra over 
$\CC ^N$.
If $\phi $ is a unital q-positive definite functional on $\cR (\cA )$, then
$\phi (X_k^2)=0$ for $k=1,\ldots ,N$. This and the q-positivity
of $\phi $ force $\phi (X_{\sigma })=0$, therefore there is only one 
q-positive functional on $\cR (\cA )$ which is not strictly q-positive.
Therefore there is no theory of orthogonal polynomials over this $\cA $.
However, the situation is different for 
$\cA =\{X_kX_l+X_lX_k \mid k,l=1,\ldots ,2N\}$. This and other 
polynomial relations will be analysed elsewhere.


\begin{thebibliography}{99}
\frenchspacing

\bibitem{AM}
J. ~Agler and J. ~E. ~McCarthy
{\em Pick interpolation and Hilbert function spaces},
Graduate Studies in Mathematics, Vol. 44, 
Amer. Math. Soc., Providence, Rhode Island, 2002.

\bibitem{ADD}
D. ~Alpay, P. ~Dewilde and H. ~Dym, 
Lossless inverse scattering and reproducing kernels for 
upper triangular operators, 
in Operator Theory: Advances and Applications,  
Vol. 47, Birkh\"auser, pp. 61-135, 1990.

\bibitem{AV}
D. ~Alpay and D. ~Volok, 
Point evaluation and Hardy space on an homogeneous tree,
lanl, OA/0309262.

\bibitem{Ar0}
W. ~Arveson, Interpolation problems in nest algebras, 
{\it J. Funct. Anal.}, {\bf 3}(1975), 208-233.

\bibitem{Ar}
W. ~Arveson, Subalgebras of $C^*$-algebras. III: Multivariable operator theory,
{\it Acta. Math.}, {\bf 181}(1998),476-514.

\bibitem{BaG}
J. ~A. ~Ball and I. ~ Gohberg, 
A commutant lifting theorem for triangular matrices with diverse 
applications, {\it Integral Equations Operator Theory},
{\bf 8}(1985),  205-267.


\bibitem{BanC}
T. ~Banks and T. ~Constantinescu,
Orthogonal polynomials in several non-commuting variables. II, 
in preparation.

\bibitem{BC}
M. ~Barakat and T. ~Constantinescu, 
Tensor Algebras and Displacement Structure. 
III. Asymptotic properties, 
lanl, FA/0306051.

\bibitem{BBCGNW}
M. ~Basseville, A. ~Benveniste, K. ~C. ~Chou, S. ~A. ~Golden, R. ~Nikoukhah,
and A. ~S. ~Willsky, 
Modeling and estimation of multiresolution stochastic processes,
{\it IEEE Trans. Inform. Theory}, {\bf 38}(1992), 766-784.

\bibitem{BBW}
M. ~Basseville, A. ~Benveniste,
and A. ~S. ~Willsky,
Multiscale autoregressive Processes, Part I: Schur-Levinson parametrizations,
{\it IEEE Trans. Signal Processing}, {\bf 40}(1992), 1915-1934;
Part II: Lattice structures for whitening and modeling,
{\it IEEE Trans. Signal Processing}, {\bf 40}(1992), 1935-1954.

\bibitem{BG}
A. ~Ben-Artzi and I. ~Gohberg, 
Orthogonal polynomials over Hilbert modules,  
in Operator Theory: Advances and Applications,  
Vol. 73, Birkh\"auser, pp. 96-126, 1994.


\bibitem{Co0}
T.~Constantinescu,
Modeling of time-variant linear systems,
INCREST preprint No.60/1985.

\bibitem{Co2}
T. ~Constantinescu,
Schur analysis of positive block matrices, 
in Operator Theory: Advances and Applications,  
Vol. 18, Birkh\"auser, pp. 191-206, 1986.

\bibitem{Co4}
T. ~Constantinescu,
Factorization of positive-definite kernels, in 
Operator Theory: Advances and Applications,  
Vol. 48, Birkh\"auser, pp. 245-260, 1990.

\bibitem{Co1}
T. ~Constantinescu,
Orthogonal polynomials in several non-commuting variables. I, 
lanl, FA/020533; to appear in I. ~Colojoara Anniversary Volume.

\bibitem{Co}
T.~Constantinescu,
{\em Schur Parameters, Factorization and Dilation Problems},
Birkh\"auser, 1996.
 
\bibitem{CJ}
T. ~Constantinescu and J. ~L. ~Johnson,
Tensor algebras and displacement structure. II. 
Non-commutative Szeg\"o theory,
{\it Zeit. f\"ur Anal. Anw.}, {\bf 21}(2002), 611--626.

\bibitem{CG}
T. ~Constantinescu and A. ~Gheondea, 
Representations of Hermitian kernels by means of Krein spaces, 
{\it Publ. RIMS}, {\bf 33}(1997), 917-951;
II. Invariant kernels, {\it Commun. Math. Phys.}, {\bf 216}(2001), 409-430.

\bibitem{DGK}
Ph. ~Delsarte, Y. ~Genin and  Y. Kamp, 
On the Toeplitz embedding of arbitrary matrices, 
{\it Linear Algebra
Appl.}, {\bf 51}(1983), 97-119.

\bibitem{Dr}
V.~Drensky, 
{\em Free Algebras and PI-Algebras}, 
Springer, 1999.

\bibitem{DX}
C. ~H. ~Dunkl and Y. ~Xu,
{\em Orthogonal Polynomials of Several Variables},
Cambridge Univ. Press, 2001.  

\bibitem{EK}
D. ~E. ~Evans and Y. ~Kawahigashi, 
{\em Quantum symmetris on operator algebras}, 
Oxford Mathematical Monographs, Clarendon Press, 1998.


\bibitem{FFGK}
C. ~Foias, A. ~E. ~Frazho, I. ~Gohberg and M. ~A. ~Kaashoek,
{\em Metric Constrained Interpolation, Commutant Lifting and Systems},
Birkh\"auser, 1998. 

\bibitem{Fr}
A. ~E. ~Frazho, 
On stochastic bilinear systems, in 
{\em Modeling and Applications of Stochastic Processes}
(U.B.Desai, Ed.), pp. 215--241, Kluwer Academic, 1988.

\bibitem{GS}
U. ~Grenander and G. ~Szeg\"o,
{\em Toeplitz Forms and their Applications}, Univ. of California Press,
California, 1958.

\bibitem{HH}
W. W. Hart and W. L. Hart, 
{\em Plane Trigonometry, Solid Geometry, and Spherical Trigonometry}, 
D. C. Heath and Co., 1942.

\bibitem{He}
P. ~Henrici, {\em  Applied and computational complex analysis.
Volume 1: Power series, integration, conformal mapping, location of zeros},
Wiley-Interscience, 1974. 

\bibitem{La}
E. ~C. ~Lance, 
{\em Hilbert $C^*$-modules}, Cambridge University Press, 1995.

\bibitem{Li}
H.~Li, 
{\em Noncommutative Gr\"obner Bases and Filtered-Graded Transfer},
LNM 1795, Springer, 2002.

\bibitem{Pa}
K.~R.~Parthasarathy,
{\em An Introduction to Quantum Stochastic Calculus},
Birkh\"auser, 1992.

\bibitem{RR}
M. ~Rosenblum and J. Rovnyak, 
{\em Hardy Classes and Operator  Theory},
Oxford Univ. Press, 1985.

\bibitem{Sz}
G. ~Szeg\"o, {\em Orthogonal Polynomials},
Colloquium Publications, 
{\bf 23}, Amer. Math. Soc., Providence, Rhode Island, 1939.

\bibitem{NF}
B. ~Sz.-Nagy and C. Foias, 
{\em Harmonic Analysis of Operators
on Hilbert Space}, North Holland, 1970.


\end{thebibliography}
\end{document}